\documentclass[11pt]{article}
\usepackage{amssymb,amsmath,booktabs,hyperref,mathtools}
\DeclareMathOperator*{\Sup}{sup}
\usepackage{graphics,authblk}
\usepackage{graphicx}
\usepackage{subcaption}
\usepackage{multirow}
\usepackage{tabularx,enumerate}
\usepackage{rotating}
\usepackage{dcolumn}
\usepackage{pdflscape}
\usepackage{array}
\usepackage{lscape}
\usepackage{anysize}
\usepackage{color}
\usepackage{amsthm}
\usepackage{listings}
\newcommand{\mathleft}{\@fleqntrue\@mathmargin0pt}
\newcommand{\mathcenter}{\@fleqnfalse}
\providecommand{\keywords}[1]{\textbf{\textit{Keywords: }} #1}
\usepackage{float}
\usepackage{soulutf8}
\usepackage[T1]{fontenc}
\usepackage{adjustbox}
\usepackage{colortbl}
\usepackage{blindtext,nameref}

\newtheorem{theorem}{Theorem}
\newtheorem{proposition}{Proposition}
\newtheorem*{definition*}{Definition}
\newtheorem{lemma}{Lemma}
\theoremstyle{definition}
\newtheorem{definition}{Definition}
\newtheorem*{ivt}{Intermediate Value Theorem (IVT)}
\newtheorem*{operation}{Arithmatic operations of fuzzy numbers using $\alpha$-cut interval}
\newtheorem*{fdc}{First decomposition theorem}
\theoremstyle{remark}
\newtheorem*{remark}{Remark}
\graphicspath{{photos/}}
\allowdisplaybreaks
\begin{document}
	\title{An $\alpha$-cut intervals based fuzzy best-Worst method for Multi-Criteria Decision-Making}
	\author{Harshit M Ratandhara, Mohit Kumar}
	\date{}
	\affil{Department of Basic Sciences,\\ Institute of Infrastructure, Technology, Research and Management, Ahmedabad, Gujarat-380026, India\\ Email: harshitratandhara1999@gmail.com, mohitiitr3@gmail.com}
	\maketitle
	\begin{abstract}
		The Best-Worst Method (BWM) is a well-known Multi-Criteria Decision-Making (MCDM) method used to calculate criteria-weights in many real-life applications. It was observed that the decision judgments used to calculate weights in BWM may be imprecise due to human involvement. To incorporate this ambiguity into the weight calculation, Guo \& Zhao proposed a model of BWM using fuzzy sets, known as Fuzzy BWM (FBWM). Although this model is known to have wide applicability, it has several limitations. One of the biggest limitations of this existing model is that the lower, modal and upper values of the fuzzy judgment are used in the weight calculation and the other values remain unused. This makes the model independent of the shape of the fuzzy number. To solve this limitation and optimize the entire shape, we propose a model of FBWM based on $\alpha$-cut intervals. This helps in reducing information loss. It turns out that although it is possible to optimize the entire shape simultaneously, it is difficult to do so due to the involvement of infinitely many constraints in the optimization problem used to calculate optimal weights. Therefore, we approximate optimal weights using finite subset, say $F$, of $[0,1]$. We then develop a technique to measure the Degree of Approximation (DoA) of a weight set and obtain a weight set with the desired DoA. For a given $F$, approximate weights are calculated using a minimization problem that has a non-linear nature and thus may lead to multiple weights. To solve this issue, we first compute the collection of all approximate weights of the criterion, which is an interval, and then adopt the center of this interval as the approximate weight of the criterion. To measure the accuracy of a weight set, we develop the concepts of Consistency Index (CI) and Consistency Ratio (CR) for the proposed model. To measure the accuracy of the weight set, we develop the concepts of consistency index (CI) and consistency ratio (CR) for the proposed model. Finally, we discuss some numerical examples and a real-world application of the proposed model in ranking of risk factors in supply chain 4.0 and compare the results with existing models.
		
	\end{abstract}
	\keywords{Multi-criteria decision-making, Fuzzy best-worst method, $\alpha$-cut intervals, Consistency index, Supply chain 4.0}	
	\section{Introduction}
	Decision making is an essential part of everyday life. To deal with decision problems, a branch of operations research called Multi-Criteria Decision-Making (MCDM) has been developed. MCDM methods help Decision Makers (DM) in decision making, especially in complex situations involving a large number of decision criteria and alternatives. AHP\cite{saaty1990make}, TOPSIS\cite{hwang1981methods}, ELECTRE\cite{roy1990outranking}, PROMETHEE\cite{brans1985note}, etc. are well-known MCDM methods.\\\\
	The Best-Worst Method (BWM) is one of the latest MCDM methods\cite{rezaei2015best}. Rezaei proposed it as a development on AHP. It is used in the calculation of weights of decision-criteria. In recent times, it has been widely used in real-world applications like supplier selection \cite{vahidi2018sustainable}, automotive\cite{van2017battle}, eco-industrial parks\cite{zhao2018comprehensive}, airline industry\cite{rezaei2017multi}, etc. It requires pairwise comparisons between criteria in calculating weights. These comparisons are provided by the DM. Because of this human involvement, there is a chance for ambiguity in these decision data. \\\\
	In 1965, Lotfi Zadeh\cite{ZADEH1965338} proposed the concept of fuzzy sets to deal with situations such as data uncertainty or information scarcity. Fuzzy sets are extension of classical sets (crisp sets). Since fuzzy sets include inherent uncertainty, they are better than classical sets in dealing with imprecise data and hence, fuzzy MCDM methods such as fuzzy AHP \cite{mikhailov2004evaluation}, fuzzy Electre \cite{sevkli2010application}, fuzzy TOPSIS \cite{chen2000extensions}, etc. have been developed. With this motivation, Guo \& Zhao\cite{guo2017fuzzy} proposed a model of BWM using fuzzy sets known as Fuzzy Best-Worst Method (FBWM). Some other models of BWM using fuzzy sets\cite{mohtashami2021novel, dong2021fuzzy}, intuitionistic fuzzy sets\cite{wan2021novel, mou2016intuitionistic}, hesitant fuzzy sets\cite{ali2019hesitant} and other extensions of classical sets are also developed. In this paper, we propose a model of FBWM using $\alpha$-cut intervals ($\alpha$-FBWM) to address some of the limitations of the existing model of FBWM\cite{guo2017fuzzy}.\\\\
	In FBWM, optimal weights are calculated using a minimization problem which is formulated using only the lower, modal and upper values of fuzzy numbers and therefore, other values with non-zero membership have no role in the weight calculation. This may result in significant loss of information. Another drawback of this approach is that the weight calculation process becomes independent of the shape of the fuzzy number. Therefore, the fuzzy numbers having different shapes but the same lower, modal and upper values will give the same weight set, which is inappropriate. In the proposed model, we optimize the entire shape of the fuzzy numbers using $\alpha$-cut intervals. This helps in reducing information loss. Another limitation of the existing model is that in weight calculation, approximate fuzzy arithmetic operations are used and, as discussed in Section 3, the difference between the exact and the approximate operations can be quite significant. Since the exact fuzzy operations are defined using $\alpha$-cut intervals, the proposed model naturally includes exact operations in weight calculation which is another advantage of the proposed model.\\\\
	In the proposed model, optimal weights are obtained using a minimization problem which is formulated by modifying the minimization problem of the existing model. This modified problem involves infinitely many constraints. Therefore, this problem is not always easy to solve, although it has an optimal solution as we prove in Subsection 4.1. To deal with this problem, we approximate optimal weights using a finite subset $F$ of $[0,1]$, i.e., we optimize a finite number of $\alpha$-cut intervals instead of optimizing all $\alpha$-cut intervals simultaneously. In particular, for $F=\{0,1\}$, an approximate weight set is precisely an optimal weight of the existing model. Then we measure the Degree of Approximation (DoA) of a weight set. A lower value of DoA indicates a better approximation. The DoA of an weight set of the existing model is one. We can obtain better weight sets than this set by considering $F$ with the higher cardinality. In fact, we can get a weight set with the desired DoA. This indicates the superiority of the proposed model than the existing model. \\\\
	In the proposed model, approximate weights are Triangular Fuzzy Numbers (TFN). These TFNs are converted to crisp weights using a defuzzification method called Graded Mean Integration Representation (GMIR). We observe that for a fix $F$, we may get multiple approximate weight sets due to the non-linear nature of the minimization problem used in calculation of weights. To deal with the same problem in crisp BWM, Rezaei\cite{rezaei2016best} developed a technique to calculate the interval-weight of a criterion and then consider the mean value of this interval as the weight of that criterion. Following a similar method, we first establish the fact that, for a given $F$, the collection of all approximate defuzzified weights of a criterion is an interval. Then after calculating this interval, the center value of this interval is adopted as the approximate weight of that criterion for $F$.\\\\
	For crisp BWM, Rezaei\cite{rezaei2015best} developed the concepts of Consistency Index (CI) and Consistency Ratio (CR) to measure the accuracy of a weight set. Guo \& Zhao\cite{guo2017fuzzy} extended these concepts to FBWM. In Section 3, we prove that this extension is not well-defined. We then modify these concepts to be well-defined in fuzzy environment. For this purpose, we derive the necessary conditions for Fuzzy Pairwise Comparison System (FPCS) to be consistent and calculate a lower bound of CI using these conditions. This lower bound leads to an upper bound of CR, which is sufficient to check the acceptability of a weight set in most cases. Finally we discuss some numerical examples and a real-world application of the proposed model in ranking risk factors in supply chain 4.0 and compare the results with the existing model. This comparison shows that the proposed model is better than the existing model as it provides better approximation of optimal weights than the existing model.\\\\
	The rest of the paper is structured as follows: In Section 2, some basic definitions and statements are discussed. Some limitations of FBWM are discussed in Section 3. In Section 4, the proposed model and some numerical examples are explained in detail. Section 5 describes a real-life application in ranking risk factors in supply chain 4.0 and Section 6 includes conclusive comments and future directions.
	\section{Preliminaries}
	\begin{definition*}\cite{ZADEH1965338}\cite{zimmermann2011fuzzy}
		Let $X$ be a universal set. Then a fuzzy set $\tilde{A}$ on $X$ is a set of ordered pairs defined as \[\tilde{A}=\{(x,\mu_{\tilde{A}}(x)): x\in X\},\] where $\mu_{\tilde{A}}:X \rightarrow [0,1]$. Here $\mu_{\tilde{A}}$ is called the membership function of $\tilde{A}$ and $\mu_{\tilde{A}}(x)$ is called the membership value of $x$ in $\tilde{A}$.
	\end{definition*}
	\begin{definition*}\cite{zimmermann2011fuzzy}
		Let $\tilde{A}$ be a fuzzy set on $X$. Then support of $\tilde{A}$, denoted as $S(\tilde{A})$, is defined as \[S(\tilde{A})=\{x\in X: \mu_{\tilde{A}}(x)>0\}.\]
	\end{definition*}
	\begin{definition*}\cite{ZADEH1965338}\cite{zimmermann2011fuzzy}
		Let $\tilde{A}$ be a fuzzy set on $X$ and $\alpha\in[0,1].$ Then \[A_\alpha=\{x\in X: \mu_{\tilde{A}}(x)\geq\alpha\}\] is called $\alpha$-cut of $\tilde{A}$.
	\end{definition*}
	\begin{remark}
		Note that if $\alpha_1\leq \alpha_2$, then $A_{\alpha_1}\supseteq A_{\alpha_2}$, i.e., $\alpha$-cuts of a fuzzy set $\tilde{A}$ forms a family of nested crisp sets.
	\end{remark}
	\begin{fdc}\cite{klir1996fuzzy}
		Let $\tilde{A}$ be a fuzzy set on $X$ . Then \[\tilde{A}=\bigcup_{\alpha\in[0,1]}{_\alpha\tilde{A}},\] 
		\begin{equation*}
			\text{where}\ \mu_{_\alpha\tilde{A}}(x)=\alpha \mu_{A_\alpha}(x)=
			\begin{cases}
				\alpha, & \text{if}\  x\in A_\alpha\\
				0, & \text{otherwise}.
			\end{cases}       
		\end{equation*}
	\end{fdc}
	\begin{remark}
		By first decomposition theorem, it follows that any fuzzy set can be determined completely with the help of its $\alpha$-cuts. Also, two fuzzy sets $\tilde{A}$ and $\tilde{B}$ on $X$ are equal iff $S(\tilde{A})=S(\tilde{B})$ and $A_\alpha=B_\alpha$ for all $\alpha\in(0,1].$
	\end{remark}
	\begin{definition*}\cite{klir1996fuzzy}
		A fuzzy set $\tilde{A}$ on $\mathbb{R}$ is said to be fuzzy number if 
		\begin{enumerate}
			\item There is $x\in X$ such that $\mu_{\tilde{A}}(x)=1$;
			\item $A_\alpha$ is closed interval for every $\alpha \in (0,1];$
			\item $S(\tilde{A})$ is bounded set.
		\end{enumerate}
	\end{definition*}
	\begin{definition*}\cite{guo2017fuzzy}                                                             
		A fuzzy number $\tilde{A}$ is said to be Triangular Fuzzy Number (TFN) if it has membership function of the form 
		\begin{equation}
			\mu_{\tilde{A}}(x)=
			\begin{cases}
				\frac{x-a}{b-a}, & \text{if} \ a\leq x\leq b \\
				\frac{c-x}{c-b}, & \text{if} \ b\leq x\leq c \\
				0, & \text{otherwise}
			\end{cases}
		\end{equation}  
		where $a,b$ and $c$ are real numbers such that $a\leq b\leq c$. A TFN is denoted as $\tilde{A}= (a,b,c).$                                     
	\end{definition*}
	\begin{operation}\cite{klir1996fuzzy}\\
		Let $\tilde{A}$ and $\tilde{B}$ be two fuzzy numbers, and let * be one of the four basic arithmetic operations: addition, substraction, multiplication and division. Define $(A*B)_\alpha=A_\alpha*B_\alpha$ for $\alpha\in(0,1].$ By decomposition theorem, \[\tilde{A}*\tilde{B}=\bigcup_{\alpha\in[0,1]}(A*B)_\alpha=\bigcup_{\alpha\in[0,1]}(A_\alpha*B_\alpha).\]
	\end{operation}
	\begin{remark}
		Since $\alpha$-cuts are crisp sets, $A_\alpha*B_\alpha$ is usual set operations between subsets of $\mathbb{R}$. Also, we need to make sure that $0\notin B_\alpha$ for $\alpha\in(0,1]$, i.e, $\mu_{\tilde{B}}(0)=0$ in case of division.
	\end{remark}
	\begin{definition*}\cite{guo2017fuzzy}
		Let $\tilde{A}=(a,b,c)$ be a TFN. Then Graded Mean Integration Representation (GMIR), denoted by $R(\tilde{A})$, is defined as
		\begin{equation}\label{GMIR}
			R(\tilde{A})=\frac{a+4b+c}{6}.
		\end{equation}
	\end{definition*} 
	\begin{ivt}\cite{royden1988real}
		Let $f$ be a real-valued continuous function on $[a,b]$ for which $f(a)<c<f(b)$. Then there is a $a<x_0<b$ such that $f(x_0)=c$.
	\end{ivt}
	\begin{remark}
		It follows from IVT that if $f$ and $g$ be real valued continuous functions on $[a,b]$ such that $f(a)>g(a)$ and $f(b)<g(b)$, then there is a $a<x_0<b$ such that $f(x_0)=g(x_0)$.
	\end{remark}
	\begin{definition*}\cite{rezaei2015best}
		Let $C=\{c_1,c_2,...,c_n\}$ be a set of decision criteria. Then fuzzy pairwise comparision matrix $\tilde{A}$ is given by 
		\begin{equation*}
			\tilde{A}=
			\begin{bmatrix}
				\tilde{a}_{11} & \tilde{a}_{12} & \cdots & \tilde{a}_{1n}\\
				\tilde{a}_{21} & \tilde{a}_{22} & \cdots & \tilde{a}_{2n}\\
				\vdots & \vdots & \ddots & \vdots \\
				\tilde{a}_{n1} & \tilde{a}_{n2} & \cdots & \tilde{a}_{nn}
			\end{bmatrix},
		\end{equation*} 
		where $\tilde{a}_{ij}$ is a fuzzy number that shows the relative preference of $i^{th}$ criterion over $j^{th}$ criterion.
	\end{definition*}
	\section{Limitations of Fuzzy Best-Worst Method}
	To deal with ambiguity in decision judgments, Guo \& Zhao\cite{guo2017fuzzy} developed a model of BWM using fuzzy sets, known as Fuzzy Best-Worst Method (FBWM). Despite the popularity of this model, it has several limitations.
	\begin{enumerate}[(i)]
		\item In FBWM, optimal weights are calculated by solving a minimization problem. The formulation of this problem depends only on the lower, modal and upper values of fuzzy comparisons, and other values with non-zero membership remain completely unused. In other words, this problem and consequently optimal weights become independent of the shape of fuzzy comparison values. 
		\item In calculation of optimal weights, They used approximate fuzzy arithmetic operations. To observe the difference between the exact value and the approximate value, consider $\tilde{A}= (-1,1,3)$ and $\tilde{B}= (1,3,5)$. So, we get approximate $\frac{\tilde{A}}{\tilde{B}}=(-1,\frac{1}{3},3)$ and exact $\frac{\tilde{A}}{\tilde{B}}(x)=\begin{cases}
			\frac{x+1}{2-2x}, \text{ if } -1\leq x\leq 0\\
			\frac{5x+1}{2x+2}, \text{ if } 0\leq x\leq \frac{1}{3}\\
			\frac{3-x}{2x+2}, \text{ if } \frac{1}{3}\leq x\leq 3\\
			0, \quad\ \ \text{ otherwise }
		\end{cases}$.\\
		These fuzzy numbers are plotted in Figure \ref{division}. Note that the difference between these values is quite significant, which directly affects the resultant weights.
		\begin{figure}[H]
			\centering
			\includegraphics[width=15cm]{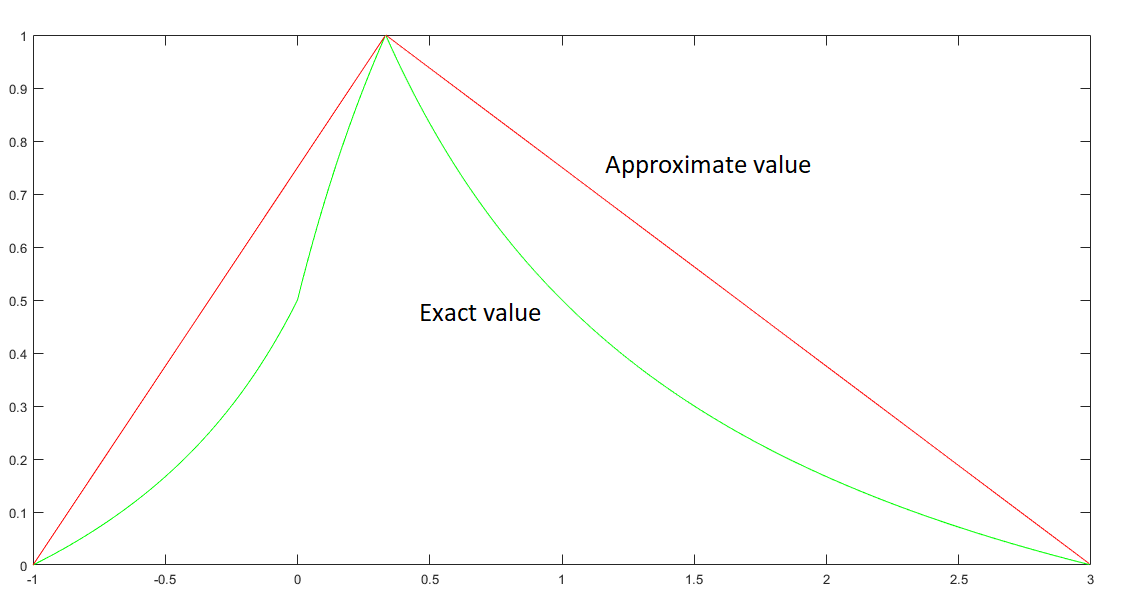}
			\caption{Exact and approximate values of $\frac{\tilde{A}}{\tilde{B}}(x)$}
			\label{division}		
		\end{figure}
		\item In crisp BWM, a consistent Pairwise Comparison System (PCS) $(A_b,A_w)$ is characterized by the condition $a_{bi}\times a_{iw}=a_{bw}$ for all $i\in\{1,2,...,n\}\setminus\{b,w\}$ because this is necessary as well as sufficient condition for the existence of accurate weights\cite{rezaei2015best}. Based on this condition, Rezaei\cite{rezaei2015best} calculated the values of the Consistency Index (CI) which are useful for measuring the accuracy of a weight set. This condition can be obtained using the equation $\frac{w_b}{w_i}\times \frac{w_i}{w_w}=\frac{w_b}{w_w}$.\\\\
		Analogously, Guo \& Zhao considered the condition $\tilde{a}_{bi}\times \tilde{a}_{iw}=\tilde{a}_{bw}$ for all $i\in\{1,2,...,n\}\setminus\{b,w\}$ as the characterization of consistent Fuzzy PCS (FPCS) $(\tilde{A}_b,\tilde{A}_w)$. But note that for a fuzzy number $\tilde{A}$, $\frac{\tilde{A}}{\tilde{A}}$ may not be $\tilde{1}$. In fact, even approximate value of $\frac{\tilde{A}}{\tilde{A}}$ may not be $\tilde{1}$. For example, consider $\tilde{A}=(1,2,3)$. Then the exact and approximate values of $\frac{\tilde{A}}{\tilde{A}}$ are shown in Figure \ref{identity}. So, the equation $\frac{\tilde{w}_b}{\tilde{w}_i}\times \frac{\tilde{w}_i}{\tilde{w}_w}=\frac{\tilde{w}_b}{\tilde{w}_w}$ is not true. Therefore, the condition $\tilde{a}_{bi}\times \tilde{a}_{iw}=\tilde{a}_{bw}$ does not ensure the existence of accurate fuzzy weights and becomes meaningless as the characterization of consistent FPCS. The values of CI for FBWM are also not well-defined as they are calculated using this condition.
		\begin{figure}[H]
			\centering
			\includegraphics[width=15cm]{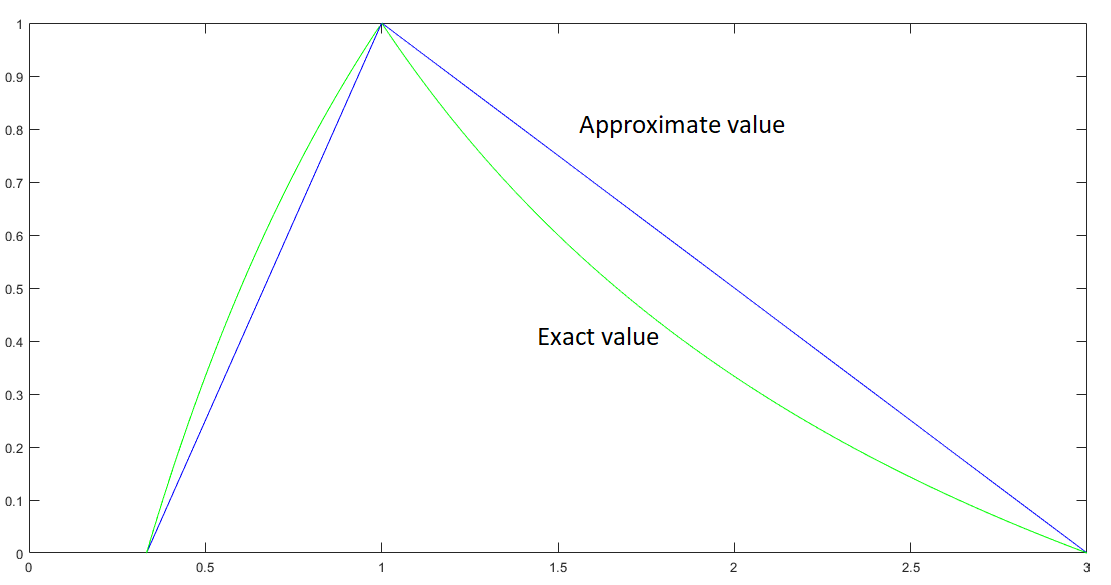}
			\caption{Exact value of $\frac{\tilde{A}}{\tilde{A}}(x)$}
			\label{identity}		
		\end{figure}
	\end{enumerate}
	\section{Fuzzy Best-Worst Method based on $\alpha$-cut intervals ($\alpha$-FBWM)}
	In this section, we develop a model of FBWM using $\alpha$-cut interval. For $\alpha\in [0,1]$, the $\alpha$-cut interval of fuzzy number $\tilde{A}$ is the collection of all elements of universal set $X$ having membership greater than or equal to $\alpha$ in $\tilde{A}$. In this model, we compare two fuzzy numbers in terms of their $\alpha$-cut intervals. This approach has two main advantages: the entire shape of fuzzy numbers, i.e., each element having non-zero membership value is included in calculation of weights and since exact fuzzy arithmetic operations are defined using $\alpha$-cut intervals, the exact values are used instead of approximate values\cite{klir1996fuzzy}. First we formulate a minimization problem to calculate an optimal weight set using $\alpha$-cut intervals. Then we derive the necessary conditions for FPCS $(\tilde{A}_b,\tilde{A}_w)$ to be consistent.
	\subsection{Calculation of optimal weights of decision criteria}
	This subsection describes the steps of $\alpha$-FBWM for calculation of optimal weights of decision criteria.\\\\
	Step 1: \textit{Formation of set of decision criteria.}\\
	In this step, Decision Maker (DM) identifies all criteria involved in decision making, i.e., decision criteria and forms their set. Let $C=\{c_1,c_2,...,c_n\}$ be a set of decision criteria.\\\\
	Step 2: \textit{Selection of the best and the worst criterion.}\\
	From $C$, DM selects the most preferable and the least preferable criterion known as the best and the worst criterion, denoted by $c_b$ and $c_w$ respectively.\\\\
	Step 3: \textit{Fuzzy pairwise comparison system}.\\
	DM provides the relative preferences of the best criterion over each criterion and relative preferences of each criterion over the worst criterion in the form of linguistic terms suggesting the degree of relative preference of one criterion over other, which are associated with TFNs using a scale given in Table \ref{tab:table1} \cite{mohtashami2021novel}. These preference values form the $b^{th}$ row $\tilde{A}_b=(\tilde{a}_{b1}, \tilde{a}_{b2}, ...,\tilde{a}_{bn})$ and $w^{th}$ column $\tilde{A}_w={(\tilde{a}_{1w}, \tilde{a}_{2w}, ..., \tilde{a}_{nw})}^T$ of fuzzy pairwise comparision matrix $\tilde{A}$. These vectors are jointly called Fuzzy Pairwise Comparison System (FPCS), denoted as $(\tilde{A}_b,\tilde{A}_w)$.
	\begin{table}[ht]
		\caption{\cite{mohtashami2021novel} Linguistic terms and associated TFNs\label{tab:table1}}
		\centering
		\begin{tabular}{|c||c|}
			\hline
			Linguistic term & TFN\\
			\hline
			\hline
			Equally preference & $\tilde{1}=(1,1,1)$\\
			\hline
			Weakly preference & $\tilde{3}=(2,3,4)$\\
			\hline
			Essentially preference & $\tilde{5}=(4,5,6)$\\
			\hline
			Very strong preference & $\tilde{7}=(6,7,8)$\\
			\hline
			Absolutely preference & $\tilde{9}=(9,9,9)$\\
			\hline
			\multirow{4}{*}{Intermediate values}& $\tilde{2}=(1,2,3)$\\
			\cline{2-2}
			& $\tilde{4}=(3,4,5)$\\
			\cline{2-2}
			& $\tilde{6}=(5,6,7)$\\
			\cline{2-2}
			& $\tilde{8}=(7,8,9)$\\
			\hline
		\end{tabular}
	\end{table}\\
	Step 4: \textit{Calculation of optimal weights}.\\
	Since $\tilde{a}_{ij}$ represents the relative preference of $i^{th}$ criterion over $j^{th}$ criterion, an optimal weight set $\{\tilde{w}_1,\tilde{w}_2,....,\tilde{w}_n\}$ is a soultion of the system of equations
	\begin{eqnarray}\label{best_to_worst}
		\frac{\tilde{w}_{b}}{\tilde{w}_{i}} = \tilde{a}_{bi},\quad \frac{\tilde{w}_{i}}{\tilde{w}_{w}} = \tilde{a}_{iw},\quad \frac{\tilde{w}_{b}}{\tilde{w}_{w}} = \tilde{a}_{bw}
	\end{eqnarray}
	for all $i\in\{1,2,...,n\}\setminus\{b,w\}$\cite{guo2017fuzzy}.\\\\
	Note that (\ref{best_to_worst}) is a system of 2n-3 fuzzy equations. This gives 
	\begin{eqnarray}\label{best_to_worst_alpha}
		\left(\frac{\tilde{w}_{b}}{\tilde{w}_{i}}\right)_\alpha = \frac{(\tilde{w}_{b})_\alpha}{(\tilde{w}_{i})_\alpha} = (\tilde{a}_{bi})_\alpha,\quad 
		\left(\frac{\tilde{w}_{i}}{\tilde{w}_{w}}\right)_\alpha = \frac{(\tilde{w}_{i})_\alpha}{(\tilde{w}_{w})_\alpha} = (\tilde{a}_{iw})_\alpha, \quad \left(\frac{\tilde{w}_{b}}{\tilde{w}_{w}}\right)_\alpha = 				\frac{(\tilde{w}_{b})_\alpha}{(\tilde{w}_{w})_\alpha} = (\tilde{a}_{bw})_\alpha
	\end{eqnarray}
	for all $i\in\{1,2,...,n\}\setminus\{b,w\}$ and $\alpha \in (0,1]$, where $\tilde{A}_\alpha$ is the $\alpha$-cut interval of $\tilde{A}$. For each $\alpha$, (\ref{best_to_worst_alpha}) gives a system of 2n-3 crisp equations. Since an $\alpha$-cut interval of a fuzzy number is closed interval, both sides of (\ref{best_to_worst_alpha}) are closed intervals.\\\\
	By property of support, we get
	\begin{eqnarray}\label{support_best_to_worst}
		S\left(\frac{\tilde{w}_{b}}{\tilde{w}_{i}}\right) = \frac{S(\tilde{w}_{b})}{S(\tilde{w}_{i})} = S(\tilde{a}_{bi}), \quad
		S\left(\frac{\tilde{w}_{i}}{\tilde{w}_{w}}\right) = \frac{S(\tilde{w}_{i})}{S(\tilde{w}_{w})} = S(\tilde{a}_{iw}),\quad
		S\left(\frac{\tilde{w}_{b}}{\tilde{w}_{w}}\right) = \frac{S(\tilde{w}_{b})}{S(\tilde{w}_{w})} = S(\tilde{a}_{bw}).
	\end{eqnarray}
	Let $(\tilde{w}_{i})_\alpha = [w_i^l(\alpha), w_i^u(\alpha)]$, $S(\tilde{w}_{i}) = [w_i^l(0), w_i^u(0)]$, $(\tilde{a}_{ij})_\alpha = [a_{ij}^l(\alpha), a_{ij}^u(\alpha)]$, $S(\tilde{a}_{ij}) = [a_{ij}^l(0), a_{ij}^u(0)]$ for all $i,j\in\{1,2,...,n\}$ and $\alpha\in(0,1]$. Then by (\ref{best_to_worst_alpha}) and (\ref{support_best_to_worst}), we get
	\begin{eqnarray}
		\left[\frac{w_b^l(\alpha)}{w_i^u(\alpha)},\frac{w_b^u(\alpha)}{w_i^l(\alpha)}\right] &=& [a_{bi}^l(\alpha),a_{bi}^u(\alpha)],\quad\ \  \left[\frac{w_b^l(0)}{w_i^u(0)},\frac{w_b^u(0)}{w_i^l(0)}\right] = [a_{bi}^l(0),a_{bi}^u(0)], \nonumber \\
		\left[\frac{w_i^l(\alpha)}{w_w^u(\alpha)},\frac{w_i^u(\alpha)}{w_w^l(\alpha)}\right] &=& [a_{iw}^l(\alpha),a_{iw}^u(\alpha)],\quad \left[\frac{w_i^l(0)}{w_w^u(0)},\frac{w_i^u(0)}{w_w^l(0)}\right] = [a_{iw}^l(0),a_{iw}^u(0)], \nonumber \\
		\left[\frac{w_b^l(\alpha)}{w_w^u(\alpha)},\frac{w_b^u(\alpha)}{w_w^l(\alpha)}\right] &=& [a_{bw}^l(\alpha),a_{bw}^u(\alpha)],\quad\left[\frac{w_b^l(0)}{w_w^u(0)},\frac{w_b^u(0)}{w_w^l(0)}\right] = [a_{bw}^l(0),a_{bw}^u(0)] \nonumber
	\end{eqnarray}
	for all $i\in\{1,2,...,n\}\setminus\{b,w\}$ and $\alpha\in(0,1].$ This gives
	\begin{eqnarray}
		\frac{w_b^l(\alpha)}{w_i^u(\alpha)} &=& a_{bi}^l(\alpha), \quad \frac{w_b^u(\alpha)}{w_i^l(\alpha)} = a_{bi}^u(\alpha), \nonumber \\
		\frac{w_i^l(\alpha)}{w_w^u(\alpha)} &=& a_{iw}^l(\alpha), \quad		\frac{w_i^u(\alpha)}{w_w^l(\alpha)} = a_{iw}^u(\alpha), \nonumber \\
		\label{equality}
		\frac{w_b^l(\alpha)}{w_w^u(\alpha)} &=& a_{bw}^l(\alpha), \quad \frac{w_b^u(\alpha)}{w_w^l(\alpha)} = a_{bw}^u(\alpha)
	\end{eqnarray}
	for all $i\in\{1,2,...,n\}\setminus\{b,w\}$ and $\alpha\in[0,1].$\\\\
	It is important to note that though $\alpha=0$ is included in (\ref{equality}), it does not represent $\alpha$-cut interval. From now on, $\alpha=0$ represents the support of a fuzzy number.
	\begin{definition}
		A weight set $\{\tilde{w}_1,\tilde{w}_2,...,\tilde{w}_n\}$ is said to be accurate if it satisfies the system of equations (\ref{equality}).
	\end{definition} 
	It may happen that (\ref{equality}) does not have a solution, i.e., accurate weight set does not exist. In such cases, we compute an optimal weight set.
	\begin{definition}
		A weight set $\{\tilde{w}_1,\tilde{w}_2,...,\tilde{w}_n\}$ having the minimum possible value of $\max\{|\frac{w_b^l(\alpha)}{w_i^u(\alpha)}- a_{bi}^l(\alpha)|, |\frac{w_b^u(\alpha)}{w_i^l(\alpha)} - a_{bi}^u(\alpha)|,|\frac{w_i^l(\alpha)}{w_w^u(\alpha)} - a_{iw}^l(\alpha)|,|\frac{w_i^u(\alpha)}{w_w^l(\alpha)} - a_{iw}^u(\alpha)|,|\frac{w_b^l(\alpha)}{w_w^u(\alpha)} - a_{bw}^l(\alpha)|,|\frac{w_b^u(\alpha)}{w_w^l(\alpha)} - a_{bw}^u(\alpha)|:i=1,2,...,n,i\neq b, i\neq w \text{ and } \alpha\in[0,1]\}$ is said to be an optimal weight set.
	\end{definition}
	Consider the following minimization problem.
	\begin{equation}\label{general_minimization}
		\begin{split}
			&\min \epsilon\\
			&\text{subject to:}\\
			&\left| \frac{w_{b}^l(\alpha)}{w_{i}^u(\alpha)} - a_{bi}^l(\alpha) \right| \leq \epsilon, \quad \left| \frac{w_{b}^u(\alpha)}{w_{i}^l(\alpha)} - a_{bi}^u(\alpha) \right| \leq \epsilon, \\
			&\left| \frac{w_{i}^l(\alpha)}{w_{w}^u(\alpha)} - a_{iw}^l(\alpha) \right| \leq \epsilon, \quad \left| \frac{w_{i}^u(\alpha)}{w_{w}^l(\alpha)} - a_{iw}^u(\alpha) \right| \leq \epsilon, \\
			&\left| \frac{w_{b}^l(\alpha)}{w_{w}^u(\alpha)} - a_{bw}^l(\alpha) \right| \leq \epsilon,\quad \left| \frac{w_{b}^u(\alpha)}{w_{w}^l(\alpha)} - a_{bw}^u(\alpha) \right| \leq \epsilon   \\
			&\text{for } i=1,2,...,n,\  i\neq b, i\neq w;\\
			&0 \leq w_{q}^l(\alpha) \leq w_{q}^u(\alpha)\quad\quad\quad\quad\quad\quad\quad\quad\ \  \text{(well-defineness of $\alpha$-cut intervals)}\\
			&\text{for } q=1,2,...,n \text{ and } \alpha\in[0,1];\\
			&N(\tilde{w}_1,\tilde{w}_2,...,\tilde{w}_n); \quad \quad \quad \quad\quad \quad \quad \quad\quad\ \text{(normalizaion condition)}\\
			&w_{r}^l(\alpha_1) \leq w_{r}^l(\alpha_2), \ w_{r}^u(\alpha_1) \geq w_{r}^u(\alpha_2)\quad \text{(nested property of $\alpha$-cut intervals)}\\
			&\text{for } r=1,...,n \text{ and } \alpha_1,\alpha_2\in[0,1] \text{ such that } \alpha_1\leq \alpha_2.
		\end{split}
	\end{equation}
	First we prove that the problem (\ref{general_minimization}) has an optimal solution. Take $\tilde{w}_i=(\frac{1}{n},\frac{1}{n},\frac{1}{n})$ for all i$\in\{1,2,...,n\}$. So, we get 
	\begin{eqnarray*}
		\left| \frac{w_{b}^l(\alpha)}{w_{i}^u(\alpha)} - a_{bi}^l(\alpha) \right|= \left| \frac{w_{b}^l(\alpha)}{w_{i}^u(\alpha)} - 1 - a_{bi}^l(\alpha) +1\right|\leq \left| \frac{w_{b}^l(\alpha)}{w_{i}^u(\alpha)} - 1\right| + \left| a_{bi}^l(\alpha) -1\right|=&0+\left| a_{bi}^l(\alpha) -1\right|\\
		&\leq \left| a_{bw}^u(\alpha) -1\right|.
	\end{eqnarray*}
	Similarly, it can be checked that, for a given weight set, each absolute difference is less than or equal to $\left| a_{bw}^u(\alpha) -1\right|$. So, we get $\max\{|\frac{w_b^l(\alpha)}{w_i^u(\alpha)}- a_{bi}^l(\alpha)|, |\frac{w_b^u(\alpha)}{w_i^l(\alpha)} - a_{bi}^u(\alpha)|,|\frac{w_i^l(\alpha)}{w_w^u(\alpha)} - a_{iw}^l(\alpha)|,|\frac{w_i^u(\alpha)}{w_w^l(\alpha)} - a_{iw}^u(\alpha)|,|\frac{w_b^l(\alpha)}{w_w^u(\alpha)} - a_{bw}^l(\alpha)|,|\frac{w_b^u(\alpha)}{w_w^l(\alpha)} - a_{bw}^u(\alpha)|:i=1,2,...,n,i\neq b, i\neq w \text{ and } \alpha\in[0,1]\}\leq \left| a_{bw}^u(\alpha) -1\right|$ and since $\epsilon^*$ is the minimum possible value of this supremum, we get $0\leq \epsilon^*\leq \left| a_{bw}^u(\alpha) -1\right|$. So, the problem (\ref{general_minimization}) has an optimal solution.\\\\
	Note that an optimal solution of the problem (\ref{general_minimization}) gives optimal $\alpha$-cut intervals for all $\alpha\in [0,1]$. Using these $\alpha$-cut intervals, an optimal weight set $\{\tilde{w}_1^*,\tilde{w}_2^*,...,\tilde{w}_n^*\}$ can be calculated using first decomposition theorem. The optimal objective value $\epsilon^*$ is a measurement of the accuracy of an optimal weight set. \\\\
	Observe that the problem (\ref{general_minimization}) has infinitely many constraints. Also it is a non-linear minimization problem and therefore may have multiple optimal solutions. This makes it difficult to solve. To deal with this issue, instead of calculating optimal weight sets, we approximate them using finite subsets of $[0,1]$. It is also difficult to compute a fuzzy number from its $\alpha$-cut intervals. Therefore, we only use triangular fuzzy weights to approximate optimal weight sets.\\\\
	Let $F=\{0=\alpha_1,\alpha_2,...,\alpha_m=1\}$ be a finite subset of $[0,1]$ such that $\alpha_1<\alpha_2<...<\alpha_m$. Let $||F||_\infty=\max\{\alpha_{i+1}-\alpha_i: i=1,2,...,m-1\}$. Here, $||F||_\infty$ is called the mesh of $F$. Now we shall optimize only those $\alpha$-cut intervals for which $\alpha\in F$ instead of all $\alpha\in[0,1]$. This will give approximate optimal weight sets corresponding to $F$.\\\\
	Let $\tilde{w}_i = (w_i^l,w_i^m,w_i^u)$, for all $i=1,2,...,n.	$ Then for $\alpha\in(0,1]$, we get $(\tilde{w}_i)_\alpha = [w_i^l+\alpha(w_i^m-w_i^l), w_i^u-\alpha(w_i^u-w_i^m)]$ and	$(\tilde{w}_i)_0 = [w_i^l,w_i^u].$ 
	By (\ref{equality}), we get
	\begin{eqnarray}
		\frac{w_b^l+\alpha(w_b^m-w_b^l)}{w_i^u-\alpha(w_i^u-w_i^m)} &=& a_{bi}^l(\alpha), \quad  \frac{w_b^u-\alpha(w_b^u-w_b^m)}{w_i^l+\alpha(w_i^m-w_i^l)} = a_{bi}^u(\alpha), \nonumber \\
		\frac{w_i^l+\alpha(w_i^m-x_i^l)}{w_w^u-\alpha(w_w^u-w_w^m)} &=& a_{iw}^l(\alpha), \quad
		\frac{w_i^u-\alpha(w_i^u-w_i^m)}{w_w^l+\alpha(w_w^m-w_w^l)} = a_{iw}^u(\alpha), \nonumber \\
		\label{equality_tfn}
		\frac{w_b^l+\alpha(w_b^m-w_b^l)}{w_w^u-\alpha(w_w^u-w_w^m)} &=& a_{bw}^l(\alpha), \quad 
		\frac{w_b^u-\alpha(w_b^u-w_b^m)}{w_w^l+\alpha(w_w^m-w_w^l)} = a_{bw}^u(\alpha) 
	\end{eqnarray}
	for all $\alpha \in [0,1]$ and $i\in\{1,2,...,n\}\setminus\{b,w\}$.\\\\
	For $\alpha \in F$ and under the assumption that the weights are TFN, the problem (\ref{general_minimization}) takes the following form.
	\begin{equation}
		\begin{split}
			\label{optimization}
			&\min \epsilon_F \\
			&\text{subject to:}\\
			&\left| \frac{w_b^l+\alpha(w_b^m-w_b^l)}{w_i^u-\alpha(w_i^u-w_i^m)} - a_{bi}^l(\alpha) \right| \leq \epsilon_F,\quad 
			\left| \frac{w_b^u-\alpha(w_b^u-w_b^m)}{w_i^l+\alpha(w_i^m-w_i^l)} - a_{bi}^u(\alpha) \right| \leq \epsilon_F, \\
			&\left| \frac{w_i^l+\alpha(w_i^m-w_i^l)}{w_w^u-\alpha(w_w^u-w_w^m)} - a_{iw}^l(\alpha) \right| \leq \epsilon_F, \quad
			\left| \frac{w_i^u-\alpha(w_i^u-w_i^m)}{w_w^l+\alpha(w_w^m-w_w^l)} - a_{iw}^u(\alpha) \right| \leq \epsilon_F, \\
			&\left| \frac{w_b^l+\alpha(w_b^m-w_b^l)}{w_w^u-\alpha(w_w^u-w_w^m)} - a_{bw}^l(\alpha) \right| \leq \epsilon_F, \quad
			\left| \frac{w_b^u-\alpha(w_b^u-w_b^m)}{w_w^l+\alpha(w_w^m-w_w^l)} - a_{bw}^u(\alpha) \right| \leq \epsilon_F, \\
			&0 \leq w_{q}^l \leq w_{q}^m \leq w_{q}^u,\quad \text{ (well-defineness and non-negetivity of weights)}\\
			&\sum_{i=1}^n R(\tilde{w}_i)=1\quad \text{ (normalization)}
		\end{split}
	\end{equation}
	for $i\in\{1,2,...,n\}\setminus\{b,w\}$, $\alpha\in F$ and $q\in\{1,2,...n\}$.\\\\
	Note that (\ref{optimization}) is a minimization problem with $3n+1$ variables and has an optimal solution of the form $(\tilde{w}_1,\tilde{w}_2,...,\tilde{w}_n,\epsilon_F^*)$. Here $\{\tilde{w}_1,\tilde{w}_2,...,\tilde{w}_n\}$ is an approximate optimal weight set corresponding to $F$, $\epsilon_F^*$ is a measurement of the accuracy of it and the set $\{R(\tilde{w}_1),R(\tilde{w}_2),...,R(\tilde{w}_n)\}$ is corresponding defuzzified approximate optimal weight set, where $R(\tilde{w})$ is GMIR of $\tilde{w}$ given by (\ref{GMIR}).\\\\ Now we estimate the Degree of Approximation (DoA) of a weight set.
	\begin{theorem}\label{doa}
		Let $(\tilde{A}_b,\tilde{A}_w)$ be FPCS formed using the scale given in Table \ref{tab:table1}, let $F=\{\alpha_1,\alpha_2,...,\alpha_m\}$ $\subset [0,1]$ such that $0=\alpha_1<\alpha_2<...<\alpha_m=1$, let $\tilde{W}_F$ be an approximate weight set corresponding to $F$, let $\eta^*=\max\{|\frac{w_b^l(\alpha)}{w_i^u(\alpha)}- a_{bi}^l(\alpha)|, |\frac{w_b^u(\alpha)}{w_i^l(\alpha)} - a_{bi}^u(\alpha)|,	|\frac{w_i^l(\alpha)}{w_w^u(\alpha)} - a_{iw}^l(\alpha)|,\frac{w_i^u(\alpha)}{w_w^l(\alpha)} - a_{iw}^u(\alpha)|,|\frac{w_b^l(\alpha)}{w_w^u(\alpha)} - a_{bw}^l(\alpha)|,|\frac{w_b^u(\alpha)}{w_w^l(\alpha)} - a_{bw}^u(\alpha)|: i=1,2,...,n,i\neq b, i\neq w, \alpha\in [0,1]\}$, and let $\epsilon^*$ be the optimal objective value of problem (\ref{general_minimization}). Then $| \epsilon^*-\eta^*|\leq||F||_\infty$, where $||F||_\infty$ is the mesh of $F$. In words, we shall express this property by saying that the Degree of Approximation (DoA) of $\tilde{W}_F$ is $||F||_\infty$.
	\end{theorem}
	\begin{proof}
		First we prove that $\eta_F^*\leq \epsilon^*$, where $\eta_F^*$ is the optimal objective value of problem (\ref{optimization}) corresponding to $F$. \\\\
		Since $\epsilon^*$ is the optimal objective value of problem (\ref{general_minimization}), there exist $\tilde{W}=\{\tilde{w}_1,\tilde{w}_2,...,\tilde{w}_n\}$ such that $\epsilon^*=\max\{|\frac{w_b^l(\alpha)}{w_i^u(\alpha)}- a_{bi}^l(\alpha)|, |\frac{w_b^u(\alpha)}{w_i^l(\alpha)} - a_{bi}^u(\alpha)|,	|\frac{w_i^l(\alpha)}{w_w^u(\alpha)} - a_{iw}^l(\alpha)|,|\frac{w_i^u(\alpha)}{w_w^l(\alpha)} - a_{iw}^u(\alpha)|,|\frac{w_b^l(\alpha)}{w_w^u(\alpha)} - a_{bw}^l(\alpha)|,|\frac{w_b^u(\alpha)}{w_w^l(\alpha)} - a_{bw}^u(\alpha)|: i=1,2,...,n,i\neq b, i\neq w, \alpha\in [0,1]\}$. This gives $\max\{|\frac{w_b^l(\alpha)}{w_i^u(\alpha)}- a_{bi}^l(\alpha)|, |\frac{w_b^u(\alpha)}{w_i^l(\alpha)} - a_{bi}^u(\alpha)|,	|\frac{w_i^l(\alpha)}{w_w^u(\alpha)} - a_{iw}^l(\alpha)|,|\frac{w_i^u(\alpha)}{w_w^l(\alpha)} - a_{iw}^u(\alpha)|,|\frac{w_b^l(\alpha)}{w_w^u(\alpha)} - a_{bw}^l(\alpha)|,|\frac{w_b^u(\alpha)}{w_w^l(\alpha)} - a_{bw}^u(\alpha)|: i=1,2,...,n,i\neq b, i\neq w, \alpha\in F\}\leq\epsilon^*$. Since $\eta_F^*$ is such smallest non-negative number, we get $\eta_F^*\leq\epsilon^*$.\\\\ 
		Let $\tilde{W}_F=\{\tilde{w}_{1F},\tilde{w}_{2F},...,\tilde{w}_{nF}\}$. Since $\tilde{W}_F$ is an approximate optimal weight set corresponding to $F$, we get  $\max\{|\frac{w_{bF}^l(\alpha)}{w_{iF}^u(\alpha)}- a_{bi}^l(\alpha)|, |\frac{w_{bF}^u(\alpha)}{w_{iF}^l(\alpha)} - a_{bi}^u(\alpha)|,	|\frac{w_{iF}^l(\alpha)}{w_{wF}^u(\alpha)} - a_{iw}^l(\alpha)|,|\frac{w_{iF}^u(\alpha)}{w_{wF}^l(\alpha)} - a_{iw}^u(\alpha)|,|\frac{w_{bF}^l(\alpha)}{w_{wF}^u(\alpha)} - a_{bw}^l(\alpha)|,|\frac{w_{bF}^u(\alpha)}{w_{wF}^l(\alpha)} - a_{bw}^u(\alpha)|: i=1,2,...,n,i\neq b, i\neq w, \alpha\in F\}=\eta_F^*\leq\eta_F^*+||F||_\infty\leq \epsilon^*+||F||_\infty$. This gives $\max\{|\frac{w_{bF}^l(\alpha)}{w_{iF}^u(\alpha)}- a_{bi}^l(\alpha)|, |\frac{w_{bF}^u(\alpha)}{w_{iF}^l(\alpha)} - a_{bi}^u(\alpha)|,	|\frac{w_{iF}^l(\alpha)}{w_{wF}^u(\alpha)} - a_{iw}^l(\alpha)|,|\frac{w_{iF}^u(\alpha)}{w_{wF}^l(\alpha)} - a_{iw}^u(\alpha)|,|\frac{w_{bF}^l(\alpha)}{w_{wF}^u(\alpha)} - a_{bw}^l(\alpha)|,|\frac{w_{bF}^u(\alpha)}{w_{wF}^l(\alpha)} - a_{bw}^u(\alpha)|: i=1,2,...,n,i\neq b, i\neq w\}\leq \epsilon^*+||F||_\infty$ for all $\alpha\in F$.\\\\
		Let $\alpha\in [0,1]\setminus F$. So, there exist $\alpha_j,\alpha_{j+1}\in F$ such that $\alpha_j<\alpha<\alpha_{j+1}$. Let $\epsilon\geq 0$ be such that $\frac{w_{bF}^l(\alpha)}{w_{iF}^u(\alpha)}=a_{bi}^l(\alpha)\pm \epsilon$. Since $\frac{w_{bF}^l(\alpha)}{w_{iF}^u(\alpha)}$ is an increasing function of $\alpha$, we get $\frac{w_{bF}^l(\alpha_j)}{w_{iF}^u(\alpha_j)}<\frac{w_{bF}^l(\alpha)}{w_{iF}^u(\alpha)}<\frac{w_{bF}^l(\alpha_{j+1})}{w_{iF}^u(\alpha_{j+1})}$. Since $\tilde{W}_F$ is an approximate optimal weight set corresponding to $F$, we get $|\frac{w_{bF}^l(\alpha)}{w_{iF}^u(\alpha)}- a_{bi}^l(\alpha)|\leq\eta_F^*\leq \epsilon^*$. In particular, we have $|\frac{w_{bF}^l(\alpha_j)}{w_{iF}^u(\alpha_j)}- a_{bi}^l(\alpha_j)|\leq \epsilon^*$ and $|\frac{w_{bF}^l(\alpha_{j+1})}{w_{iF}^u(\alpha_{j+1})}- a_{bi}^l(\alpha_{j+1})|\leq \epsilon^*$. This implies $a_{bi}^l(\alpha_j)-\epsilon^*\leq\frac{w_{bF}^l(\alpha_j)}{w_{iF}^u(\alpha_j)}$ and $\frac{w_{bF}^l(\alpha_{j+1})}{w_{iF}^u(\alpha_{j+1})}\leq a_{bi}^l(\alpha_{j+1})+\epsilon^*$. So, $a_{bi}^l(\alpha_j)-\epsilon^*<a_{bi}^l(\alpha)\pm \epsilon<a_{bi}^l(\alpha_{j+1})+\epsilon^*$. Now observe that for the scale given in Table \ref{tab:table1}, $a_{bi}^l(\alpha_1)-a_{bi}^l(\alpha_2)\leq \alpha_1-\alpha_2$ for all $\alpha_1,\alpha_2\in [0,1]$ such $\alpha_1<\alpha_2$. We also have $\alpha_{j+1}-\alpha_j\leq ||F||_\infty$ for all $\alpha_j,\alpha_{j+1}\in F$. This gives $-||F||_\infty -\epsilon^*\leq \frac{w_{bF}^l(\alpha)}{w_{iF}^u(\alpha)}- a_{bi}^l(\alpha)\leq ||F||_\infty+\epsilon^*$ and thus $|\frac{w_{bF}^l(\alpha)}{w_{iF}^u(\alpha)}- a_{bi}^l(\alpha)|\leq \epsilon^*+||F||_\infty$. Similarly, it can be proven that each absolute difference is less than or equal to $\epsilon^*+||F||_\infty$ for all $\alpha\in [0,1]\setminus F$. This gives $\max\{|\frac{w_{bF}^l(\alpha)}{w_{iF}^u(\alpha)}- a_{bi}^l(\alpha)|, |\frac{w_{bF}^u(\alpha)}{w_{iF}^l(\alpha)} - a_{bi}^u(\alpha)|,	|\frac{w_{iF}^l(\alpha)}{w_{wF}^u(\alpha)} - a_{iw}^l(\alpha)|,|\frac{w_{iF}^u(\alpha)}{w_{wF}^l(\alpha)} - a_{iw}^u(\alpha)|,|\frac{w_{bF}^l(\alpha)}{w_{wF}^u(\alpha)} - a_{bw}^l(\alpha)|,|\frac{w_{bF}^u(\alpha)}{w_{wF}^l(\alpha)} - a_{bw}^u(\alpha)|: i=1,2,...,n,i\neq b, i\neq w\}\leq \epsilon^*+||F||_\infty$ for all $\alpha\in [0,1]\setminus F$.\\\\
		From the above discussion, we get $0\leq \eta^*\leq \epsilon^*+||F||_\infty $, i.e., $| \epsilon^*-\eta^*|\leq||F||_\infty$. Hence the proof.
	\end{proof}
	\begin{remark}
		Note that for $F=\{0,1\}$, problem (\ref{optimization}) is the minimization problem of the existing model of FBWM\cite{guo2017fuzzy} and therefore, $\tilde{W}_F$ is an optimal solution of the existing model. Now, Theorem \ref{doa} implies that this weight set has a DoA equal to $1$. This implies that an optimal weight set of the existing model is one of the approximate weight sets of the proposed model with DoA equal to $1$, which is very high. We can get a better approximation by choosing $F$ with a lower value of $||F||_\infty$. For example, $F=\{0,0.25,0.5,0.75,1\}$ gives a weight set with DoA equal to $0.25$, which is better than an optimal weight set of the existing model. This indicates the superiority of the proposed model over the existing model.
	\end{remark}
	\subsection{Interval-weights}
	Rezaei\cite{rezaei2015best} proposed the non-linear model of BWM for calculation of optimal weights of decision criteria. One of the limitations of this model is that it may lead to multiple optimal weight sets. To deal with this problem, Rezaei\cite{rezaei2016best} developed two minimization problems that give the GLB and the LUB of the collection of optimal weights of a criterion, which is an interval and then the average value of interval-weight is considered as the weight of that criterion. In $\alpha$-FBWM also, minimization problem (\ref{optimization}) is non-linear, which may lead to multiple solutions. So, in this subsection, we extend this approach for the proposed model. First we establish the fact that for the given $F$, the collection of all defuzzified approximate optimal weights of a criterion is an interval. If it is unique, say $w$, then it represents the interval $[w,w]$. Otherwise the following theorem implies that it is an interval.
	\begin{theorem}\label{Interval}
		Let $C=\{c_1,c_2,...,c_n\}$ be a set of decision criteria with $c_b$ and $c_w$ as the best and the worst criterion respectively, let $(\tilde{A}_b,\tilde{A}_w)$ be a FPCS, let $F$ be a finite subset of $[0,1]$, and let $w_{i_0}'$ and $w_{i_0}''$ be defuzzified approximate optimal weights of criterion $c_{i_0}$  corresponding to $F$ for some $i_0\in\{1,2,...,n\}$ such that $w_{i_0}'<w_{i_0}''$. Then each $w_{i_0}\in[w_{i_0}',w_{i_0}'']$ is an defuzzified approximate optimal weight of $c_{i_0}$ corresponding to $F$.
	\end{theorem}
	\begin{proof}
		Since $w_{i_0}'$ and $w_{i_0}''$ are defuzzified approximate optimal weights of $c_{i_0}$ corresponding to $F$, there exist $\tilde{W}'=\{\tilde{w}_1',\tilde{w}_2',...,\tilde{w}_n'\}$ and $\tilde{W}''=\{\tilde{w}_1'',\tilde{w}_2'',...,\tilde{w}_n''\}$ such that $R(\tilde{w}_{i_0}')=w_{i_0}'$, $R(\tilde{w}_{i_0}'')=w_{i_0}''$,
		\begin{eqnarray*}
			&&\left| \frac{{w_b'}^l+\alpha({w_b'}^m-{w_b'}^l)}{{w_i'}^u-\alpha({w_i'}^u-{w_i'}^m)} - a_{bi}^l(\alpha) \right| \leq \epsilon^*,\quad 
			\left| \frac{{w_b'}^u-\alpha({w_b'}^u-{w_b'}^m)}{{w_i'}^l+\alpha({w_i'}^m-{w_i'}^l)} - a_{bi}^u(\alpha) \right| \leq \epsilon^*,  \nonumber\\
			&&\left| \frac{{w_i'}^l+\alpha({w_i'}^m-{w_i'}^l)}{{w_w'}^u-\alpha({w_w'}^u-{w_w'}^m)} - a_{iw}^l(\alpha) \right| \leq \epsilon^*, \quad
			\left| \frac{{w_i'}^u-\alpha({w_i'}^u-{w_i'}^m)}{{w_w'}^l+\alpha({w_w'}^m-{w_w'}^l)} - a_{iw}^u(\alpha) \right| \leq \epsilon^*, \nonumber\\
			&&\left| \frac{{w_b'}^l+\alpha({w_b'}^m-{w_b'}^l)}{{w_w'}^u-\alpha({w_w'}^u-{w_w'}^m)} - a_{bw}^l(\alpha) \right| \leq \epsilon^*, \quad
			\left| \frac{{w_b'}^u-\alpha({w_b'}^u-{w_b'}^m)}{{w_w'}^l+\alpha({w_w'}^m-{w_w'}^l)} - a_{bw}^u(\alpha) \right| \leq \epsilon^*,  \nonumber\\
			&&\left| \frac{{w_b''}^l+\alpha({w_b''}^m-{w_b''}^l)}{{w_i''}^u-\alpha({w_i''}^u-{w_i''}^m)} - a_{bi}^l(\alpha) \right| \leq \epsilon^*,\quad 
			\left| \frac{{w_b''}^u-\alpha({w_b''}^u-{w_b''}^m)}{{w_i''}^l+\alpha({w_i''}^m-{w_i''}^l)} - a_{bi}^u(\alpha) \right| \leq \epsilon^*,  \nonumber\\
			&&\left| \frac{{w_i''}^l+\alpha({w_i''}^m-{w_i''}^l)}{{w_w''}^u-\alpha({w_w''}^u-{w_w''}^m)} - a_{iw}^l(\alpha) \right| \leq \epsilon^*, \quad
			\left| \frac{{w_i''}^u-\alpha({w_i''}^u-{w_i''}^m)}{{w_w''}^l+\alpha({w_w''}^m-{w_w''}^l)} - a_{iw}^u(\alpha) \right| \leq \epsilon^*, \nonumber\\
			&&\left| \frac{{w_b''}^l+\alpha({w_b''}^m-{w_b''}^l)}{{w_w''}^u-\alpha({w_w''}^u-{w_w''}^m)} - a_{bw}^l(\alpha) \right| \leq \epsilon^*, \quad
			\left| \frac{{w_b''}^u-\alpha({w_b''}^u-{w_b''}^m)}{{w_w''}^l+\alpha({w_w''}^m-{w_w''}^l)} - a_{bw}^u(\alpha) \right| \leq \epsilon^*,  \nonumber\\
			&&0 \leq {w_q'}^l \leq {w_q'}^m \leq {w_q'}^u,\quad 0 \leq {w_q''}^l \leq {w_q''}^m \leq {w_q''}^u,\\
			&&\sum_{i=1}^{n}R(\tilde{w}_i')=1, \quad \sum_{i=1}^{n}R(\tilde{w}_i'')=1
		\end{eqnarray*}
		for $i\in\{1,2,...,n\}\setminus\{b,w\}$, $\alpha\in F$ and $q\in\{1,2,...n\}$, where $\tilde{w}_k'=({w_k'}^l,{w_k'}^m,{w_k'}^u)$, $\tilde{w}_k''=({w_k''}^l,{w_k''}^m,{w_k''}^u)$ for all $k\in\{1,2,...,n\}$ and $\epsilon^*$ is the optimal objective value of problem (\ref{optimization}).\\\\
		Since $w_{i_0}\in[w_{i_0}',w_{i_0}'']$, there exists $\lambda\in[0,1]$ such that $\lambda w_{i_0}'+(1-\lambda)w_{i_0}''=w_{i_0}$. Consider $\tilde{w}_i=\lambda \tilde{w}_i'+(1-\lambda) \tilde{w}_i''=(\lambda{w_i'}^l+(1-\lambda){w_i''}^l,\lambda{w_i'}^m+(1-\lambda){w_i''}^m,\lambda{w_i'}^u+(1-\lambda){w_i''}^u)$ for all $i\in\{1,2,...,n\}$. So, we get ${w_i}^l=\lambda{w_i'}^l+(1-\lambda){w_i''}^l$, ${w_i}^m=\lambda{w_i'}^m+(1-\lambda){w_i''}^m$, and ${w_i}^u=\lambda{w_i'}^u+(1-\lambda){w_i''}^u$. Note that $0\leq w_i^l\leq w_i^m \leq w_i^u$. So, each $\tilde{w}_i$ is a well-defined non-negative TFN. Also note that $R(\tilde{w}_{i_0})=w_{i_0}$.\\\\
		Consider $\tilde{W}=\{\tilde{w}_i:i=1,2,...,n\}$. Let $\alpha\in F$. Now
		\begin{eqnarray*}
			\frac{{w_b}^l+\alpha({w_b}^m-{w_b}^l)}{{w_i}^u-\alpha({w_i}^u-{w_i}^m)}&=&\frac{\lambda({w_b'}^l+\alpha({w_b'}^m-{w_b'}^l))+(1-\lambda)({w_b''}^l+\alpha({w_b''}^m-{w_b''}^l))}{\lambda({w_i'}^u-\alpha({w_i'}^u-{w_i'}^m))+(1-\lambda)({w_i''}^u-\alpha({w_i''}^u-{w_i''}^m))}.
		\end{eqnarray*}
		Let $f(\lambda)=\frac{\lambda({w_b'}^l+\alpha({w_b'}^m-{w_b'}^l))+(1-\lambda)({w_b''}^l+\alpha({w_b''}^m-{w_b''}^l))}{\lambda({w_i'}^u-\alpha({w_i'}^u-{w_i'}^m))+(1-\lambda)({w_i''}^u-\alpha({w_i''}^u-{w_i''}^m))}$, $\lambda\in[0,1]$.\\\\
		So, $f'(\lambda)=\frac{({w_b''}^l+\alpha({w_b''}^m-{w_b''}^l))({w_i'}^u-\alpha({w_i'}^u-{w_i'}^m))-({w_b'}^l+\alpha({w_b'}^m-{w_b'}^l))({w_i''}^u-\alpha({w_i''}^u-{w_i''}^m))}{(\lambda({w_i'}^u-\alpha({w_i'}^u-{w_i'}^m))+(1-\lambda)({w_i''}^u-\alpha({w_i''}^u-{w_i''}^m)))^2}$.\\\\
		Note that $f'(\lambda)\geq0$, i.e, $f$ is increasing in $[0,1]$ iff $\frac{{w_b''}^l+\alpha({w_b''}^m-{w_b''}^l)}{{w_i''}^u-\alpha({w_i''}^u-{w_i''}^m))}\geq\frac{{w_b'}^l+\alpha({w_b'}^m-{w_b'}^l)}{{w_i'}^u-\alpha({w_i'}^u-{w_i'}^m)}$.\\\\
		In this case, we get $\frac{{w_b'}^l+\alpha({w_b'}^m-{w_b'}^l)}{{w_i'}^u-\alpha({w_i'}^u-{w_i'}^m)}\leq\frac{{w_b}^l+\alpha({w_b}^m-{w_b}^l)}{{w_i}^u-\alpha({w_i}^u-{w_i}^m)}\leq\frac{{w_b''}^l+\alpha({w_b''}^m-{w_b''}^l)}{{w_i''}^u-\alpha({w_i''}^u-{w_i''}^m)}$.\\\\
		This implies $\frac{{w_b'}^l+\alpha({w_b'}^m-{w_b'}^l)}{{w_i'}^u-\alpha({w_i'}^u-{w_i'}^m)}-a_{bi}^l(\alpha)\leq\frac{{w_b}^l+\alpha({w_b}^m-{w_b}^l)}{{w_i}^u-\alpha({w_i}^u-{w_i}^m)}-a_{bi}^l(\alpha)\leq\frac{{w_b''}^l+\alpha({w_b''}^m-{w_b''}^l)}{{w_i''}^u-\alpha({w_i''}^u-{w_i''}^m)}-a_{bi}^l(\alpha)$.\\\\
		If $\frac{{w_b}^l+\alpha({w_b}^m-{w_b}^l)}{{w_i}^u-\alpha({w_i}^u-{w_i'}^m)}-a_{bi}^l(\alpha)\leq 0$, then $\frac{{w_b'}^l+\alpha({w_b'}^m-{w_b'}^l)}{{w_i'}^u-\alpha({w_i'}^u-{w_i'}^m)}-a_{bi}^l(\alpha)\leq0$. Since $|\cdot|$ is a decreasing function in $(-\infty,0]$, we get $$\bigg|\frac{{w_b}^l+\alpha({w_b}^m-{w_b}^l)}{{w_i}^u-\alpha({w_i}^u-{w_i}^m)}-a_{bi}^l(\alpha)\bigg|\leq\bigg|\frac{{w_b'}^l+\alpha({w_b'}^m-{w_b'}^l)}{{w_i'}^u-\alpha({w_i'}^u-{w_i'}^m)}-a_{bi}^l(\alpha)\bigg|\leq \epsilon^*.$$
		If $\frac{{w_b}^l+\alpha({w_b}^m-{w_b}^l)}{{w_i}^u-\alpha({w_i}^u-{w_i}^m)}-a_{bi}^l(\alpha)> 0$, then $\frac{{w_b''}^l+\alpha({w_b''}^m-{w_b''}^l)}{{w_i''}^u-\alpha({w_i''}^u-{w_i''}^m)}-a_{bi}^l(\alpha)>0$. Since $|\cdot|$ is an increasing function in $(0,\infty)$, we get $$\bigg|\frac{{w_b}^l+\alpha({w_b}^m-{w_b}^l)}{{w_i}^u-\alpha({w_i}^u-{w_i}^m)}-a_{bi}^l(\alpha)\bigg|\leq\bigg|\frac{{w_b''}^l+\alpha({w_b''}^m-{w_b''}^l)}{{w_i''}^u-\alpha({w_i''}^u-{w_i''}^m)}-a_{bi}^l(\alpha)\bigg|\leq \epsilon^*.$$
		Similar argument can be given if $f$ is decreasing in $[0,1]$, i.e., $\frac{{w_b''}^l+\alpha({w_b''}^m-{w_b''}^l)}{{w_i''}^u-\alpha({w_i''}^u-{w_i''}^m))}\leq\frac{{w_b'}^l+\alpha({w_b'}^m-{w_b'}^l)}{{w_i'}^u-\alpha({w_i'}^u-{w_i'}^m)}$.\\\\
		In similar manner, it can be proven that for $\tilde{W}$, all other abosulte differences are also less than or equal to $\epsilon^*$. Since $\alpha$ is arbitrary, this holds for all values of $\alpha\in F$. Hence $\tilde{W}$ is an optimal weight set. Now it remains to prove that $\tilde{W}$ is normalized.
		\begin{eqnarray*}
			\text{Now, } \sum_{i=1}^n R(\tilde{w}_i)&=& \sum_{i=1}^n \frac{(\lambda {w_i'}^l+ (1-\lambda) {w_i''}^l)+4*(\lambda {w_i'}^m+ (1-\lambda) {w_i''}^m)+(\lambda {w_i'}^u+ (1-\lambda) {w_i''}^u)}{6}\\
			&=&\frac{\lambda \sum_{i=1}^n ({w_i'}^l+4*{w_i'}^m+{w_i'}^u)+ (1-\lambda) \sum_{i=1}^n ({w_i''}^l+4*{w_i''}^m+{w_i''}^u) }{6}\\
			&=&\lambda \frac{\sum_{i=1}^n ({w_i'}^l+4*{w_i'}^m+{w_i'}^u)}{6}+ (1-\lambda) \frac{({w_i''}^l+4*{w_i''}^m+{w_i''}^u) }{6}\\
			&=&\lambda + 1-\lambda\\
			&=&1.	
		\end{eqnarray*}
		So, $\tilde{W}$ is a normalized approximate optimal weight set corresponding to $F$ with $\tilde{R}(\tilde{w}_{i_0})=w_{i_0}$. This completes the proof.
	\end{proof}
	Theorem \ref{Interval} shows that for the given $F$, the collection of defuzzified approximate optimal weights of criterion is of the form $(a,b),(a,b],[a,b)$ or $[a,b]$, i.e., an interval. For computational purposes, this interval can be treated as $[a,b]$. To calculate the values of $a$ and $b$, i.e., the GLB and the LUB of this interval, consider the following minimization problems.
	\begin{equation}	\label{optimization_glb}
		\begin{split}
			&\min R(\tilde{w}_k)\bigg(=\frac{w_k^l+4w_k^m+w_k^u}{6}\bigg) \\
			&\text{subject to:}\\
			&\left| \frac{w_b^l+\alpha(w_b^m-w_b^l)}{w_i^u-\alpha(w_i^u-w_i^m)} - a_{bi}^l(\alpha) \right| \leq \epsilon_F^*,\quad 
			\left| \frac{w_b^u-\alpha(w_b^u-w_b^m)}{w_i^l+\alpha(w_i^m-w_i^l)} - a_{bi}^u(\alpha) \right| \leq \epsilon_F^*,  \\
			&\left| \frac{w_i^l+\alpha(w_i^m-w_i^l)}{w_w^u-\alpha(w_w^u-w_w^m)} - a_{iw}^l(\alpha) \right| \leq \epsilon_F^*, \quad
			\left| \frac{w_i^u-\alpha(w_i^u-w_i^m)}{w_w^l+\alpha(w_w^m-w_w^l)} - a_{iw}^u(\alpha) \right| \leq \epsilon_F^*, \\
			&\left| \frac{w_b^l+\alpha(w_b^m-w_b^l)}{w_w^u-\alpha(w_w^u-w_w^m)} - a_{bw}^l(\alpha) \right| \leq \epsilon_F^*, \quad
			\left| \frac{w_b^u-\alpha(w_b^u-w_b^m)}{w_w^l+\alpha(w_w^m-w_w^l)} - a_{bw}^u(\alpha) \right| \leq \epsilon_F^*, \\
			&0 \leq w_{q}^l \leq w_{q}^m \leq w_{q}^u,\\
			&\sum_{i=1}^{n}R(\tilde{w}_i)= 1\\
			&\text{for } i\in\{1,2,...,n\}\setminus\{b,w\}, \alpha\in F \text{ and } q\in\{1,2,...n\} \text{ and }
		\end{split}
	\end{equation}
	\begin{equation}
		\begin{split}
			\label{optimization_lub}
			&\max R(\tilde{w}_k)\bigg(=\frac{w_k^l+4w_k^m+w_k^u}{6}\bigg) \\
			&\text{subject to:}\\
			&\left| \frac{w_b^l+\alpha(w_b^m-w_b^l)}{w_i^u-\alpha(w_i^u-w_i^m)} - a_{bi}^l(\alpha) \right| \leq \epsilon_F^*,\quad 
			\left| \frac{w_b^u-\alpha(w_b^u-w_b^m)}{w_i^l+\alpha(w_i^m-w_i^l)} - a_{bi}^u(\alpha) \right| \leq \epsilon_F^*,  \\
			&\left| \frac{w_i^l+\alpha(w_i^m-w_i^l)}{w_w^u-\alpha(w_w^u-w_w^m)} - a_{iw}^l(\alpha) \right| \leq \epsilon_F^*, \quad
			\left| \frac{w_i^u-\alpha(w_i^u-w_i^m)}{w_w^l+\alpha(w_w^m-w_w^l)} - a_{iw}^u(\alpha) \right| \leq \epsilon_F^*, \\
			&\left| \frac{w_b^l+\alpha(w_b^m-w_b^l)}{w_w^u-\alpha(w_w^u-w_w^m)} - a_{bw}^l(\alpha) \right| \leq \epsilon_F^*, \quad
			\left| \frac{w_b^u-\alpha(w_b^u-w_b^m)}{w_w^l+\alpha(w_w^m-w_w^l)} - a_{bw}^u(\alpha) \right| \leq \epsilon_F^*, \\
			&0 \leq w_{q}^l \leq w_{q}^m \leq w_{q}^u,\\
			&\sum_{i=1}^{n}R(\tilde{w}_i)= 1\\
			&\text{for } i\in\{1,2,...,n\}\setminus\{b,w\}, \alpha\in F \text{ and } q\in\{1,2,...n\},
		\end{split}
	\end{equation}
	where $\epsilon_F^*$ is the optimal objective value corresponding to $F$ in problem (\ref{optimization}).\\\\
	Note that (\ref{optimization_glb}) and (\ref{optimization_lub}) are minimization problems having $3n$ variables which have the GLB and the LUB of the collection of defuzzified approximate optimal weights of criterion $k$ corresponding to $F$ as the optimal objective values respectively. Let $w_k^*(\text{lower})$ and $w_k^*(\text{upper})$ be the said GLB and LUB. Then $w_k^*(\text{avg})=\frac{w_k^*(\text{lower})+w_k^*(\text{upper})}{2}$. Further computation can be done by treating $w_k^*(\text{avg})$ as the weight of criterion $k$.
	\subsection{Consistency index}
	To measure the accuracy of weight set, Rezaei\cite{rezaei2015best} developed the concepts of consistency index and consistency ratio for crisp BWM. Here we extend these concepts for fuzzy environment.\\\\
	For the given $\tilde{a}_{bw}$, Consistency Index (CI) is the maximum possible optimal objective value in the problem (\ref{general_minimization}). So,
	$$\text{CI}|_{\tilde{a}_{bw}}=\sup\{\epsilon^*: \epsilon^* \text{ is the optimal objective value for some FPCS } (\tilde{A}_b,\tilde{A}_w) \text{ having given }\tilde{a}_{bw} \}.$$
	Let $\{\tilde{w}_1,\tilde{w}_2,...,\tilde{w}_n\}$ be a weight set. Then Consistency Ratio (CR) for this weight set is defined as $$\text{CR}=\frac{\eta^*}{\text{CI}},$$ where $\eta^*=\max\{|\frac{w_b^l(\alpha)}{w_i^u(\alpha)}- a_{bi}^l(\alpha)|, |\frac{w_b^u(\alpha)}{w_i^l(\alpha)} - a_{bi}^u(\alpha)|,|\frac{w_i^l(\alpha)}{w_w^u(\alpha)} - a_{iw}^l(\alpha)|,|\frac{w_i^u(\alpha)}{w_w^l(\alpha)} - a_{iw}^u(\alpha)|,|\frac{w_b^l(\alpha)}{w_w^u(\alpha)} - a_{bw}^l(\alpha)|,|\frac{w_b^u(\alpha)}{w_w^l(\alpha)} - a_{bw}^u(\alpha)|:i=1,2,...,n,i\neq b, i\neq w \text{ and } \alpha\in[0,1]\}$.
	\begin{remark}
		Let $(\tilde{A}_b,\tilde{A}_w)$ be a FPCS, and let $\{\tilde{w}^*_1,\tilde{w}^*_2,...,\tilde{w}^*_n\}$ be an optimal weight set for $(\tilde{A}_b,\tilde{A}_w)$, i.e., a solution of problem (\ref{general_minimization}). Then $\epsilon^*=\eta^*$, where $\epsilon^*$ is the optimal objective value. So, $\text{CR}=\frac{\epsilon^*}{\text{CI}}.$
	\end{remark}
	\begin{definition}
		A FPCS $(\tilde{A}_b,\tilde{A}_w)$ is said to be consistent if there exist fuzzy numbers (not necessarily triangular) $\tilde{w}_1,\tilde{w}_2,...,\tilde{w}_n$ satisfying the system of equations (\ref{best_to_worst}).
	\end{definition}
	\begin{theorem}\label{consistency}
		Let $(\tilde{A}_b,\tilde{A}_w)$ be a consistent FPCS. Then for all  $i,i_1,i_2\in\{1,2,...,n\}\setminus\{b,w\},\ \alpha\in[0,1]$ and some functions $K_1,K_2$ of $\alpha$,
		\begin{eqnarray}
			\label{consistency_lower}
			&&a_{bi}^l(\alpha) \times a_{iw}^u(\alpha) = K_1(\alpha),\\
			\label{consistency_upper}
			&&a_{bi}^u(\alpha) \times a_{iw}^l(\alpha) = K_2(\alpha),\\ 
			\label{consistency_lower_upper}	
			&&a_{bi_1}^l(\alpha) \times a_{i_1w}^u(\alpha) \times a_{bi_2}^u(\alpha) \times a_{i_2w}^l(\alpha) =  a_{bw}^l(\alpha) \times a_{bw}^u(\alpha),\\
			\label{consistency_lower_upper_inequality}
			&&a_{bi}^l(\alpha) \times a_{iw}^u(\alpha) \leq a_{bw}^u(\alpha), \\
			\label{consistency_upper_lower_inequality}
			&&a_{bi}^u(\alpha) \times a_{iw}^l(\alpha) \leq a_{bw}^u(\alpha), \\
			\label{consistency_lower_lower_inequality}
			&&a_{bi}^l(\alpha) \times a_{iw}^l(\alpha) \leq a_{bw}^l(\alpha)
		\end{eqnarray}
		and the functions $f,g,h$ of $\alpha$ given by
		\begin{eqnarray}
			\label{consistency_increasing_1}
			&&f(\alpha)= \frac{a_{bi}^l(\alpha) \times a_{iw}^u(\alpha)}{a_{bw}^u(\alpha)},\\
			\label{consistency_increasing_2}
			&&g(\alpha)= \frac{a_{bi}^u(\alpha) \times a_{iw}^l(\alpha)}{a_{bw}^u(\alpha)},\\
			\label{consistency_increasing_3}
			&&h(\alpha)= \frac{a_{bi}^l(\alpha) \times a_{iw}^l(\alpha)}{a_{bw}^l(\alpha)}
		\end{eqnarray} 	
		are increasing functions.
	\end{theorem}
	\begin{proof}
		Let $(\tilde{A}_b,\tilde{A}_w)$ be a consistent FPCS, i.e., system of equations (\ref{best_to_worst}) has a solution. So, there exist $\{\tilde{w}_1,\tilde{w}_2,....,\tilde{w}_n\}$ such that
		\begin{eqnarray*}
			\frac{\tilde{w}_{b}}{\tilde{w}_{i}} = \tilde{a}_{bi},\quad \frac{\tilde{w}_{i}}{\tilde{w}_{w}} = \tilde{a}_{iw}\quad\text{and} \quad\frac{\tilde{w}_{b}}{\tilde{w}_{w}} = \tilde{a}_{bw}
		\end{eqnarray*}
		for all $i\in\{1,2,...,n\}\setminus\{b,w\}$. From (\ref{equality}), it follows that
		\begin{eqnarray*}
			&&a_{bi}^l(\alpha) \times a_{iw}^u(\alpha) =\frac{w_b^l(\alpha)}{w_w^l(\alpha)},\quad  a_{bi}^u(\alpha) \times a_{iw}^l(\alpha) =\frac{w_b^u(\alpha)}{w_w^u(\alpha)}\quad  \text{ and }\\ 
			&&a_{bi_1}^l(\alpha) \times a_{i_1w}^u(\alpha) \times a_{bi_2}^u(\alpha) \times a_{i_2w}^l(\alpha) \nonumber= a_{bw}^l(\alpha) \times a_{bw}^u(\alpha)  		
		\end{eqnarray*}
		for all $i,i_1,i_2\in\{1,2,...,n\}\setminus\{b,w\},\alpha\in[0,1]$.\\\\		
		Observe that $\frac{w_b^l(\alpha)}{w_w^l(\alpha)}$ and  $\frac{w_b^u(\alpha)}{w_w^u(\alpha)}$ are independent from $i$. Take $K_1(\alpha)=\frac{w_b^l(\alpha)}{w_w^l(\alpha)}$ and  $K_2(\alpha)=\frac{w_b^u(\alpha)}{w_w^u(\alpha)}$. This gives
		\begin{eqnarray*}
			a_{bi}^l(\alpha) \times a_{iw}^u(\alpha) =K_1(\alpha)\quad  \text{ and }\quad	a_{bi}^u(\alpha) \times a_{iw}^l(\alpha) =K_2(\alpha). 
		\end{eqnarray*}
		for all $i\in\{1,2,...,n\}\setminus\{b,w\},\alpha\in[0,1]$.\\\\	
		Now, from (\ref{equality}), we get
		\begin{eqnarray*}
			w_b^l(\alpha) &=& \frac{a_{bi}^l(\alpha) \times a_{iw}^u(\alpha)}{a_{bw}^u(\alpha)}\times w_b^u(\alpha), \\
			w_w^l(\alpha) &=& \frac{a_{bi}^u(\alpha) \times a_{iw}^l(\alpha)}{a_{bw}^u(\alpha)}\times w_w^u(\alpha), \\
			w_i^l(\alpha) &=& \frac{a_{bi}^l(\alpha) \times a_{iw}^l(\alpha)}{a_{bw}^l(\alpha)}\times  w_i^u(\alpha) 		
		\end{eqnarray*}
		for all $i\in\{1,2,...,n\}\setminus\{b,w\}$, $\alpha\in[0,1]$. We know that $w_i^l(\alpha)\leq w_i^u(\alpha) $ for all $\alpha\in[0,1]$ (well-defineness of $\alpha$-cut intervals). So, we get
		\begin{eqnarray*}
			&&a_{bi}^l(\alpha) \times a_{iw}^u(\alpha) \leq a_{bw}^u(\alpha), \\
			&&a_{bi}^u(\alpha) \times a_{iw}^l(\alpha) \leq a_{bw}^u(\alpha), \\
			&&a_{bi}^l(\alpha) \times a_{iw}^l(\alpha) \leq a_{bw}^l(\alpha)
		\end{eqnarray*}
		for all $i\in\{1,2,...,n\}\setminus\{b,w\}$ and $\alpha\in[0,1]$.\\\\
		Let $\alpha_1\leq\alpha_2$. Then $w_i^l(\alpha_1)\leq w_i^l(\alpha_2)$ and $w_i^u(\alpha_1)\geq w_i^u(\alpha_2)$ (nested property of $\alpha$-cut intevals). This gives 
		\begin{eqnarray*}
			&&\frac{a_{bi}^l(\alpha_1) \times a_{iw}^u(\alpha_1)}{a_{bw}^u(\alpha_1)}\leq \frac{a_{bi}^l(\alpha_2) \times a_{iw}^u(\alpha_2)}{a_{bw}^u(\alpha_2)},\\
			&&\frac{a_{bi}^u(\alpha_1) \times a_{iw}^l(\alpha_1)}{a_{bw}^u(\alpha_1)}\leq \frac{a_{bi}^u(\alpha_2) \times a_{iw}^l(\alpha_2)}{a_{bw}^u(\alpha_2)}\quad \text{and} \\
			&&\frac{a_{bi}^l(\alpha_1) \times a_{iw}^l(\alpha_1)}{a_{bw}^l(\alpha_1)}\leq \frac{a_{bi}^l(\alpha_2) \times a_{iw}^l(\alpha_2)}{a_{bw}^l(\alpha_2)}.
		\end{eqnarray*}
		So, the functions $f,g$ and $h$ of $\alpha$ given by
		\begin{eqnarray*}
			f(\alpha)= \frac{a_{bi}^l(\alpha) \times a_{iw}^u(\alpha)}{a_{bw}^u(\alpha)},\quad
			g(\alpha)= \frac{a_{bi}^u(\alpha) \times a_{iw}^l(\alpha)}{a_{bw}^u(\alpha)} \quad\text{and}\quad
			h(\alpha)= \frac{a_{bi}^l(\alpha) \times a_{iw}^l(\alpha)}{a_{bw}^l(\alpha)}
		\end{eqnarray*} 	
		are increasing functions for all $i\in\{1,2,...,n\}\setminus\{b,w\}$.
	\end{proof}
	Let atleast one of (\ref{consistency_lower}) to (\ref{consistency_increasing_3}) does not hold. Then $(\tilde{A}_b,\tilde{A}_w)$ be an inconsistent FPCS \\\\
	\textbf{Case 1.} Suppose that (\ref{consistency_lower}) does not hold. So, there exist $i_1,i_2\in\{1,2,...,n\}\setminus\{b,w\}$ and $\alpha\in [0,1]$ such that
	$a_{bi_1}^l(\alpha) \times a_{i_1w}^u(\alpha) \neq a_{bi_2}^l(\alpha) \times a_{i_2w}^u(\alpha)$. In order to make $(\tilde{A}_b,\tilde{A}_w)$ consistent, we need to change one or more of these four values.\\\\
	The Consistency Value (CV) of $a_{bi_1}^l(\alpha) \times a_{i_1w}^u(\alpha) \neq a_{bi_2}^l(\alpha) \times a_{i_2w}^u(\alpha)$, say $\epsilon$, is defined as 
	\begin{multline*}
		\epsilon=\inf\{\max\{\epsilon_1,\epsilon_2,\epsilon_3,\epsilon_4\}: (a_{bi_1}^l(\alpha)\pm \epsilon_1 )\times (a_{i_1w}^u(\alpha)\pm\epsilon_2) = (a_{bi_2}^l(\alpha)\pm\epsilon_3) \times (a_{i_2w}^u(\alpha)\pm\epsilon_4)\\ \text{ and }\epsilon_1,\epsilon_2,\epsilon_3,\epsilon_4\geq 0\}.
	\end{multline*} So, the CV of $a_{bi_1}^l(\alpha) \times a_{i_1w}^u(\alpha) \neq a_{bi_2}^l(\alpha) \times a_{i_2w}^u(\alpha)$ is the minimum possible change required in comparision values to make them consistent. Without loss of generality, we may assume that $a_{bi_1}^l(\alpha) \times a_{i_1w}^u(\alpha) < a_{bi_2}^l(\alpha) \times a_{i_2w}^u(\alpha)$. By IVT, it can be proven that $\epsilon$ is the root of the equation
	\begin{eqnarray*}
		(a_{bi_1}^l(\alpha) + x) \times (a_{i_1w}^u(\alpha) + x)=(a_{bi_2}^l(\alpha) - x) \times (a_{i_2w}^u(\alpha) - x).
	\end{eqnarray*} 
	This gives
	\begin{eqnarray}\label{cv_lower}
		\epsilon = \frac{a_{bi_2}^l(\alpha) \times a_{i_2w}^u(\alpha)-a_{bi_1}^l(\alpha) \times a_{i_1w}^u(\alpha)}{a_{bi_1}^l(\alpha) + a_{i_1w}^u(\alpha) + a_{bi_2}^l(\alpha) + a_{i_2w}^u(\alpha)}.
	\end{eqnarray}
	Observe that
	\begin{eqnarray*}
		\epsilon < \frac{a_{bi_2}^l(\alpha) \times a_{i_2w}^u(\alpha)}{a_{bi_1}^l(\alpha) + a_{i_1w}^u(\alpha) + a_{bi_2}^l(\alpha) + a_{i_2w}^u(\alpha)} < \frac{a_{bi_2}^l(\alpha) \times a_{i_2w}^u(\alpha)}{a_{i_2w}^u(\alpha)} = a_{bi_2}^l(\alpha), \text{ i.e., }\epsilon<a_{bi_2}^l(\alpha).
	\end{eqnarray*}
	Similarly, $\epsilon < a_{i_2w}^u(\alpha).$
	\begin{remark}
		CV can be defined even if $a_{bi_1}^l(\alpha) \times a_{i_1w}^u(\alpha) = a_{bi_2}^l(\alpha) \times a_{i_2w}^u(\alpha)$. In this case, CV $=0$.
	\end{remark}
	\begin{proposition}\label{lower_bound_1}
		Let $i_1,i_2\in\{1,2,...,n\}\setminus\{b,w\}$ and $\alpha\in [0,1]$ such that
		$a_{bi_1}^l(\alpha) \times a_{i_1w}^u(\alpha) < a_{bi_2}^l(\alpha) \times a_{i_2w}^u(\alpha)$, let $\epsilon$ be the CV of this inequality, and let $\epsilon^*$ be the optimal objective value in (\ref{general_minimization}). Then $\epsilon\leq\epsilon^*$.
	\end{proposition}
	\begin{proof}
		Since $\epsilon$ is the CV of $a_{bi_1}^l(\alpha) \times a_{i_1w}^u(\alpha) < a_{bi_2}^l(\alpha) \times a_{i_2w}^u(\alpha)$, we have $(a_{bi_1}^l(\alpha) + \epsilon) \times (a_{i_1w}^u(\alpha) + \epsilon)=(a_{bi_2}^l(\alpha) - \epsilon) \times (a_{i_2w}^u(\alpha) - \epsilon)$. As $\epsilon^*$ is the solution of (\ref{general_minimization}), there exist fuzzy numbers $\tilde{w}^*_1,\tilde{w}^*_2,...,\tilde{w}^*_n$ such that 
		\begin{eqnarray*}
			&&\left| \frac{w_{b}^{*l}(\alpha)}{w_{i}^{*u}(\alpha)} - a_{bi}^l(\alpha) \right| \leq \epsilon^*, \quad \left| \frac{w_{b}^{*u}(\alpha)}{w_{i}^{*l}(\alpha)} - a_{bi}^u(\alpha) \right| \leq \epsilon^*,\quad \left| \frac{w_{i}^{*l}(\alpha)}{w_{w}^{*u}(\alpha)} - a_{iw}^l(\alpha) \right| \leq \epsilon^*, \\
			&&\left| \frac{w_{i}^{*u}(\alpha)}{w_{w}^{*l}(\alpha)} - a_{iw}^u(\alpha) \right| \leq \epsilon^*,\quad \left| \frac{w_{b}^{*l}(\alpha)}{w_{w}^{*u}(\alpha)} - a_{bw}^l(\alpha) \right| \leq \epsilon^*,\quad \left| \frac{w_{b}^{*u}(\alpha)}{w_{w}^{*l}(\alpha)} - a_{bw}^u(\alpha) \right| \leq \epsilon^*
		\end{eqnarray*}
		for all $i\in\{1,2,...,n\}\setminus\{b,w\}$ and $\alpha\in[0,1]$. So,
		\begin{eqnarray*}
			&&\left| \frac{w_{b}^{*l}(\alpha)}{w_{i_1}^{*u}(\alpha)} - a_{bi_1}^l(\alpha) \right| =\eta_1, \quad \left| \frac{w_{i_1}^{*u}(\alpha)}{w_{w}^{*l}(\alpha)} - a_{i_1w}^u(\alpha) \right| =\eta_2,\\	&&\left| \frac{w_{b}^{*l}(\alpha)}{w_{i_2}^{*u}(\alpha)} - a_{bi_2}^l(\alpha) \right| =\eta_3, \quad \left| \frac{w_{i_2}^{*u}(\alpha)}{w_{w}^{*l}(\alpha)} - a_{i_2w}^u(\alpha) \right| =\eta_4,
		\end{eqnarray*}
		for some $0\leq\eta_1,\eta_2,\eta_3,\eta_4\leq\epsilon^*$. Thus
		\begin{eqnarray*}
			&&\frac{w_{b}^{*l}(\alpha)}{w_{i_1}^{*u}(\alpha)} = a_{bi_1}^l(\alpha) \pm \eta_1, \quad \frac{w_{i_1}^{*u}(\alpha)}{w_{w}^{*l}(\alpha)} = a_{i_1w}^u(\alpha) \pm \eta_2,\\	
			&&\frac{w_{b}^{*l}(\alpha)}{w_{i_2}^{*u}(\alpha)} = a_{bi_2}^l(\alpha) \pm\eta_3, \quad \frac{w_{i_2}^{*u}(\alpha)}{w_{w}^{*l}(\alpha)} = a_{i_2w}^u(\alpha) \pm\eta_4.
		\end{eqnarray*} This gives $(a_{bi_1}^l(\alpha) \pm \eta_1)\times(a_{i_1w}^u(\alpha) \pm \eta_2)=(a_{bi_2}^l(\alpha) \pm\eta_3)\times(a_{i_2w}^u(\alpha) \pm\eta_4)$. To prove $\epsilon\leq\epsilon^*$, it suffices to prove that at least one of $\epsilon\leq\eta_1$, $\epsilon\leq\eta_2$, $\epsilon\leq\eta_3$ and $\epsilon\leq\eta_4$ hold. Suppose, if possible, $\epsilon>\eta_1$, $\epsilon>\eta_2$, $\epsilon>\eta_3$ and $\epsilon>\eta_4$. This implies $(a_{bi_1}^l(\alpha) + \epsilon) \times (a_{i_1w}^u(\alpha) + \epsilon)>(a_{bi_1}^l(\alpha) \pm \eta_1)\times(a_{i_1w}^u(\alpha) \pm \eta_2)$ and $(a_{bi_2}^l(\alpha) - \epsilon) \times (a_{i_2w}^u(\alpha) - \epsilon)<(a_{bi_2}^l(\alpha) \pm\eta_3)\times(a_{i_2w}^u(\alpha) \pm\eta_4)$. So, $(a_{bi_2}^l(\alpha) - \epsilon) \times (a_{i_2w}^u(\alpha) - \epsilon)<(a_{bi_1}^l(\alpha) + \epsilon) \times (a_{i_1w}^u(\alpha) + \epsilon)$, which is contradiction. Hence $\epsilon\leq\epsilon^*$.
	\end{proof}	
	\begin{proposition}	\label{consistency_lower_higher}
		Let $i_1,i_2,i_3\in\{1,2,...,n\}\setminus\{b,w\}$ and $\alpha\in[0,1]$ be such that
		\begin{eqnarray*}
			a_{bi_1}^l(\alpha) \leq a_{bi_2}^l(\alpha) \leq a_{bi_3}^l(\alpha), \quad a_{i_1w}^u(\alpha) \leq a_{i_2w}^u(\alpha) \leq a_{i_3w}^u(\alpha), 
		\end{eqnarray*}
		and let $\epsilon_1,\epsilon_2,\epsilon_3$ be the CVs of
		\begin{eqnarray*}
			a_{bi_1}^l(\alpha) \times a_{i_1w}^u(\alpha) &\leq& a_{bi_2}^l(\alpha) \times a_{i_2w}^u(\alpha),\\
			a_{bi_2}^l(\alpha) \times a_{i_2w}^u(\alpha) &\leq& a_{bi_3}^l(\alpha) \times a_{i_3w}^u(\alpha),\\
			a_{bi_1}^l(\alpha) \times a_{i_1w}^u(\alpha) &\leq& a_{bi_3}^l(\alpha) \times a_{i_3w}^u(\alpha)			  
		\end{eqnarray*}
		respectively. Then $\epsilon_1 \leq \epsilon_3$ and $\epsilon_2 \leq \epsilon_3$.
	\end{proposition} 
	\begin{proof} 
		Suppose, if possible, $\epsilon_3 < \epsilon_1$. By hypothesis, we have
		$a_{bi_2}^l(\alpha) \leq a_{bi_3}^l(\alpha)$. So, $(a_{bi_2}^l(\alpha) - \epsilon_1) \leq (a_{bi_3}^l(\alpha) - \epsilon_1) $.
		As $\epsilon_3 < \epsilon_1$,
		$(a_{bi_3}^l(\alpha) - \epsilon_1) < (a_{bi_3}^l(\alpha) - \epsilon_3)$.
		This gives
		$ 0 < (a_{bi_2}^l(\alpha) - \epsilon_1) < (a_{bi_3}^l(\alpha) - \epsilon_3)$.
		Similarly,
		$ 0 < (a_{i_2w}^u(\alpha) - \epsilon_1) < (a_{i_3w}^u(\alpha) - \epsilon_3)$.
		So, we get
		$ (a_{bi_1}^l(\alpha) + \epsilon_1) \times (a_{i_1w}^u(\alpha) + \epsilon_1) =(a_{bi_2}^l(\alpha) - \epsilon_1) \times (a_{i_2w}^u(\alpha) - \epsilon_1) <(a_{bi_3}^l(\alpha) - \epsilon_3) \times (a_{i_3w}^u(\alpha) - \epsilon_3)=(a_{bi_1}^l(\alpha) + \epsilon_3) \times (a_{i_1w}^u(\alpha) + \epsilon_3)$. Thus
		$(a_{bi_1}^l(\alpha) + \epsilon_1) \times (a_{i_1w}^u(\alpha) + \epsilon_1) < (a_{bi_1}^l(\alpha) + \epsilon_3) \times (a_{i_1w}^u(\alpha) + \epsilon_3)$. But $\epsilon_3 < \epsilon_1$ implies
		$(a_{bi_1}^l(\alpha) + \epsilon_3) \times (a_{i_1w}^u(\alpha) + \epsilon_3) <  (a_{bi_1}^l(\alpha) + \epsilon_1) \times (a_{i_1w}^u(\alpha) + \epsilon_1)$,
		which is contradiction. So, $\epsilon_1 \leq \epsilon_3$. Similarly, it can be proven that $\epsilon_2 \leq \epsilon_3$. 
	\end{proof}
	Note that $\tilde{a}_{bi},\tilde{a}_{iw}\in\{\tilde{1},\tilde{2},...,\tilde{a}_{bw}\}$. So, by Propositions \ref{consistency_lower_higher} , it follows that, in this case, we get the maximum possible CV if $\tilde{a}_{bi_1}=\tilde{a}_{i_1w}=\tilde{1}$ and $\tilde{a}_{bi_2}=\tilde{a}_{i_2w}=\tilde{a}_{bw}$ for some $i_1,i_2\in\{1,2,...,n\}\setminus\{b,w\}$.	
	\begin{proposition}	\label{max_cor_value_1}
		Let $i_1,i_2\in\{1,2,...,n\}\setminus\{b,w\}$ be such that $\tilde{a}_{bi_1}=\tilde{a}_{i_1w}=\tilde{1}$ and $\tilde{a}_{bi_2}=\tilde{a}_{i_2w}=\tilde{a}_{bw}$, and let $\epsilon_\alpha$ be the CV of $1^l(\alpha) \times 1^u(\alpha) < a_{bw}^l(\alpha) \times a_{bw}^u(\alpha)$. Then $\displaystyle\Sup_{\alpha\in[0,1]}\{\epsilon_\alpha\} = \epsilon_1$. Consequently, $\epsilon_1$ is the maximum possible CV for Case 1. Also $\epsilon_1\leq$ CI, i.e, $\epsilon_1$ is a lower bound of CI.
	\end{proposition}
	\begin{proof}
		If $\tilde{a}_{bw}=\tilde{9}$, then by equation (\ref{cv_lower}), we get $\epsilon_\alpha=4$ for all $\alpha\in[0,1]$.\\\\		
		Now assume that $\tilde{a}_{bw}\neq\tilde{9}$. Let $\alpha\in[0,1]$. Observe that $a_{bw}^l(\alpha)= a_{bw}^m-1+\alpha,\ a_{bw}^u(\alpha)= a_{bw}^m+1-\alpha \ \text{and}\ \ 1^l(\alpha)=1^u(\alpha)=1$. So, by equation (\ref{cv_lower}), we get 
		\begin{eqnarray}
			\epsilon_\alpha &=& \frac{(a_{bw}^m-1+\alpha) \times (a_{bw}^m+1-\alpha)-(1 \times 1)}{1 + 1 + a_{bw}^m + a_{bw}^m} \nonumber \\
			&=& \frac{a^2_{bw}-(1-\alpha)^2-(1 \times 1)}{1 + 1 + a_{bw}^m + a_{bw}^m} \nonumber \\
			&\leq& \frac{a^2_{bw}-(1 \times 1)}{1 + 1 + a_{bw}^m + a_{bw}^m} \nonumber \\
			&=& \epsilon_1.\nonumber
		\end{eqnarray}
		So, in both cases, we get $\displaystyle\Sup_{\alpha\in[0,1]}\{\epsilon_\alpha\} = \epsilon_1$.	Now, Proposition \ref{consistency_lower_higher} and Proposition \ref{lower_bound_1} implies that $\epsilon_1$ is the maximum possible CV in Case 1 and $\epsilon_1\leq$ CI respectively. This completes the proof.
	\end{proof}
	The values of such $\epsilon_1$ corresponding to different values of $\tilde{a}_{bw}$ are given in Table \ref{consistency_indices} .\\\\
	\textbf{Case 2.} Suppose that (\ref{consistency_upper}) does not hold. So, there exist $i_1,i_2 \in\{1,2,...,n\}\setminus\{b,w\}$ and $\alpha\in[0,1]$ such that
	$a_{bi_1}^u(\alpha) \times a_{i_1w}^l(\alpha) \neq a_{bi_2}^u(\alpha) \times a_{i_2w}^l(\alpha)$. Then the CV of $a_{bi_1}^u(\alpha) \times a_{i_1w}^l(\alpha) \neq a_{bi_2}^u(\alpha) \times a_{i_2w}^l(\alpha)$, say $\epsilon$, is defined as 
	\begin{multline*}
		\epsilon=\inf\{\max\{\epsilon_1,\epsilon_2,\epsilon_3,\epsilon_4\}: (a_{bi_1}^u(\alpha)\pm \epsilon_1 )\times (a_{i_1w}^l(\alpha)\pm\epsilon_2) = (a_{bi_2}^u(\alpha)\pm\epsilon_3) \times (a_{i_2w}^l(\alpha)\pm\epsilon_4)\\ \text{ and }\epsilon_1,\epsilon_2,\epsilon_3,\epsilon_4\geq 0\}.
	\end{multline*} Without loss of generality, we may assume that $a_{bi_1}^u(\alpha) \times a_{i_1w}^l(\alpha) < a_{bi_2}^u(\alpha) \times a_{i_2w}^l(\alpha)$. By replicating the argument of Case 1, we can prove that
	\begin{eqnarray}
		\label{cv_upper}
		\epsilon &=& \frac{a_{bi_2}^u(\alpha) \times a_{i_2w}^l(\alpha)-a_{bi_1}^u(\alpha) \times a_{i_1w}^l(\alpha)}{a_{bi_1}^u(\alpha) + a_{i_1w}^l(\alpha) + a_{bi_2}^u(\alpha) + a_{i_2w}^l(\alpha)}.
	\end{eqnarray}
	\begin{proposition}\label{lower_bound_2}
		Let $i_1,i_2\in\{1,2,...,n\}\setminus\{b,w\}$ and $\alpha\in [0,1]$ such that
		$a_{bi_1}^u(\alpha) \times a_{i_1w}^l(\alpha) < a_{bi_2}^u(\alpha) \times a_{i_2w}^l(\alpha)$, let $\epsilon$ be the CV of this inequality, and let $\epsilon^*$ be the optimal objective value in (\ref{general_minimization}). Then $\epsilon\leq\epsilon^*$.
	\end{proposition}
	Similar to Case 1, it can be proven that CV of $1^u(\alpha) \times 1^l(\alpha) \neq a_{bw}^u(\alpha) \times a_{bw}^l(\alpha)$ at $\alpha=1$ is the maximum possible CV in Case 2, which is same as Case 1. These values are given in Table \ref{consistency_indices} .\\\\		
	\textbf{Case 3.} Suppose that (\ref{consistency_lower_upper}) does not hold. So, there exist $i_1,i_2 \in\{1,2,...,n\}\setminus\{b,w\}$ such that
	$a_{bi_1}^l(\alpha) \times a_{i_1w}^u(\alpha) \times a_{bi_2}^u(\alpha) \times a_{i_2w}^l(\alpha) \neq a_{bw}^l(\alpha) \times a_{bw}^u(\alpha)$ for some $\alpha\in [0,1]$. There are two possibilities.\\\\	
	\textbf{Subcase 1.} $a_{bi_1}^l(\alpha) \times a_{i_1w}^u(\alpha) \times a_{bi_2}^u(\alpha) \times a_{i_2w}^l(\alpha) > a_{bw}^l(\alpha) \times a_{bw}^u(\alpha)$. Here, the CV of $a_{bi_1}^l(\alpha) \times a_{i_1w}^u(\alpha) \times a_{bi_2}^u(\alpha) \times a_{i_2w}^l(\alpha) > a_{bw}^l(\alpha) \times a_{bw}^u(\alpha)$, say $\epsilon$, is defined as 
	\begin{multline*}
		\epsilon=\inf\{\max\{\epsilon_1,\epsilon_2,\epsilon_3,\epsilon_4,\epsilon_5,\epsilon_6\}: (a_{bi_1}^l(\alpha)\pm \epsilon_1 )\times (a_{i_1w}^u(\alpha)\pm\epsilon_2)\times(a_{bi_2}^u(\alpha)\pm\epsilon_3) \times (a_{i_2w}^l(\alpha)\pm\epsilon_4)=\\(a_{bw}^l(\alpha)\pm\epsilon_5) \times (a_{bw}^u(\alpha)\pm\epsilon_6) \text{ and }\epsilon_1,\epsilon_2,\epsilon_3,\epsilon_4,\epsilon_5,\epsilon_6\geq 0\}.
	\end{multline*}
	Consider the equation
	\begin{eqnarray}
		\label{cor_equation}
		(a_{bi_1}^l(\alpha)-x) \times (a_{i_1w}^u(\alpha)-x)\times (a_{bi_2}^u(\alpha) -x) \times (a_{i_2w}^l(\alpha)-x)=(a_{bw}^l(\alpha)+x) \times (a_{bw}^u(\alpha)+x).
	\end{eqnarray}
	Let $f(x)=(a_{bi_1}^l(\alpha)-x) \times (a_{i_1w}^u(\alpha)-x)\times (a_{bi_2}^u(\alpha) -x) \times (a_{i_2w}^l(\alpha)-x)$ and $g(x)= (a_{bw}^l(\alpha)+x) \times (a_{bw}^u(\alpha)+x)$, $x\in\mathbb{R}.$ Then $f$ and $g$ are continuous functions. Observe that $f(0)>g(0)$. Let $y=\min\{a_{bi_1}^l(\alpha),a_{i_1w}^u(\alpha),a_{bi_2}^u(\alpha),a_{i_2w}^l(\alpha)\}$. Then $f(y)<g(y)$. So, by IVT, there exist $0<c<y$ such that $f(c)=g(c)$. So, $c$ is a positive root of above equation, $c<a_{bi_1}^l(\alpha)$, $c<a_{i_1w}^u(\alpha)$, $c<a_{bi_2}^u(\alpha)$ and $c<a_{i_2w}^l(\alpha)$.\\\\
	Note that $\epsilon$ is the smallest positive root of equation (\ref{cor_equation}). So, it follows from above discussion that $\epsilon$ is strictly less than $a_{bi_1}^l(\alpha)$, $a_{i_1w}^u(\alpha)$, $a_{bi_2}^u(\alpha)$, and $a_{i_2w}^l(\alpha)$.
	\begin{proposition}\label{lower_bound_3_1}
		Let $i_1,i_2\in\{1,2,...,n\}\setminus\{b,w\}$ and $\alpha\in [0,1]$ such that
		$a_{bi_1}^l(\alpha) \times a_{i_1w}^u(\alpha) \times a_{bi_2}^u(\alpha) \times a_{i_2w}^l(\alpha) > a_{bw}^l(\alpha) \times a_{bw}^u(\alpha)$, let $\epsilon$ be the CV of this inequality, and let $\epsilon^*$ be the optimal objective value in (\ref{general_minimization}). Then $\epsilon\leq\epsilon^*$.
	\end{proposition}
	\begin{proposition}
		\label{consistency_lower_upper_lower_upper_inequality}
		Let $i_1,i_2,i_3,i_4\in\{1,2,...,n\}\setminus\{b,w\}$ be such that
		\begin{eqnarray}
			a_{bi_1}^l(\alpha) \leq a_{bi_2}^l(\alpha) , \quad a_{i_1w}^u(\alpha) \leq a_{i_2w}^u(\alpha),\quad a_{bi_3}^l(\alpha) \leq a_{bi_4}^l(\alpha) , \quad a_{i_3w}^u(\alpha) \leq a_{i_4w}^u(\alpha),  \nonumber
		\end{eqnarray}
		and let $\epsilon_1,\epsilon_2,\epsilon_3$ be the CVs of
		\begin{eqnarray*}
			a_{bi_1}^l(\alpha) \times a_{i_1w}^u(\alpha) \times a_{bi_3}^u(\alpha) \times a_{i_3w}^l(\alpha) &\geq& a_{bw}^l(\alpha) \times a_{bw}^u(\alpha), \\ 
			a_{bi_2}^l(\alpha) \times a_{i_2w}^u(\alpha) \times a_{bi_3}^u(\alpha) \times a_{i_3w}^l(\alpha) &\geq& a_{bw}^l(\alpha) \times a_{bw}^u(\alpha),\\
			a_{bi_1}^l(\alpha) \times a_{i_1w}^u(\alpha) \times a_{bi_4}^u(\alpha) \times a_{i_4w}^l(\alpha) &\geq& a_{bw}^l(\alpha) \times a_{bw}^u(\alpha)  
		\end{eqnarray*}
		respectively. Then $\epsilon_1 \leq \epsilon_2$ and $\epsilon_1 \leq \epsilon_3$.
	\end{proposition}
	In Proposition \ref{consistency_lower_upper_lower_upper_inequality}, $i_1$ and $i_2$ are not necessarily distinct. So, in this case, CV is maximum for the scale \ref{tab:table1} if $\tilde{a}_{bi}=\tilde{a}_{iw}=\tilde{a}_{bw}$ for some $i\in\{1,2,...,n\}\setminus\{b,w\}$.
	\begin{lemma}\label{upper_bound}
		Let $i\in\{1,2,...,n\}\setminus\{b,w\}$ be such that $\tilde{a}_{bi}=\tilde{a}_{iw}=\tilde{a}_{bw}$, and let $\epsilon_\alpha$ be the CV of $ a_{bi}^l(\alpha) \times a_{iw}^u(\alpha) \times a_{bi}^u(\alpha) \times a_{iw}^l(\alpha)> a_{bw}^l(\alpha) \times a_{bw}^u(\alpha)$ for $\alpha\in[0,1]$. Then $\epsilon_\alpha<a_{bw}^m-\sqrt{2}$.
	\end{lemma}
	\begin{proof}
		If $\tilde{a}_{bw}=\tilde{9}$, then $\epsilon_\alpha=5.2279<9-\sqrt{2}$ for all $\alpha\in[0,1]$ and we are done.\\\\
		Now assume that $\tilde{a}_{bw}\neq\tilde{9}$. For $\alpha\in[0,1]$, we get $a_{bw}^l(\alpha)= a_{bw}^m-1+\alpha$ and $a_{bw}^u(\alpha)= a_{bw}^m+1-\alpha$. This along with equation (\ref{cor_equation}) gives
		\begin{eqnarray*}
			((a_{bw}^m-1+\alpha)-\epsilon_\alpha)^2 \times ((a_{bw}^m+1-\alpha)-\epsilon_\alpha)^2 =((a_{bw}^m-1+\alpha)+\epsilon_\alpha) \times ((a_{bw}^m+1-\alpha)+\epsilon_\alpha).
		\end{eqnarray*}
		Fix $\alpha\in[0,1]$. Let $f(x)=((a_{bw}^m-1+\alpha)-x)^2 \times (a_{bw}^m+1-\alpha-x)^2$ and $g(x)=(a_{bw}^m-1+\alpha+x) \times ((a_{bw}^m+1-\alpha)+x)$, $x\in\mathbb{R}.$ Then $f$ and $g$ are continuous functions. Observe that $f(0)>g(0)$. Let $y=a_{bw}^m-\sqrt{2}$. Then $f(y)=(2-(1-\alpha)^2)$ and $g(y)=(2a_{bw}^m-\sqrt{2})^2 -(1-\alpha)^2$. If we show that $f(y)<g(y)$, then by IVT, there exist $0<c<y$ such that $f(c)=g(c)$. So, $c$ is a positive root of above equation such that  $c<(a_{bw}^m-\sqrt{2})$. As $\epsilon_\alpha$ is the smallest positive root, $\epsilon_\alpha<(a_{bw}^m-\sqrt{2})$. So, it remains to prove $(2-(1-\alpha)^2)<(2a_{bw}^m-\sqrt{2})^2 -(1-\alpha)^2$ for $\alpha\in[0,1].$\\\\
		Define $h_1(\alpha)=\epsilon_\alpha<(a_{bw}^m-\sqrt{2})$ and $h_2(\alpha)=(2-(1-\alpha)^2)<(2a_{bw}^m-\sqrt{2})^2 -(1-\alpha)^2$ for $\alpha\in[0,1]$. Since $a_{bw}^m\geq 2, h_2(0)\geq 5.70$. Thus $4=h_1(1)<5.70\leq h_2(0)$. Also, $h_1$ and $h_2$ are continuous, increasing functions. So, for $\alpha\in[0,1]$, $h_1(\alpha)\leq h_1(1) < h_2(0) \leq h_2(\alpha)$. This completes the proof.
	\end{proof}
	\begin{lemma}\label{dif_fun}
		Let $f(\epsilon,\alpha)= (a_{bw}^m-1+\alpha-\epsilon)^2 \times (a_{bw}^m+1-\alpha-\epsilon)^2 - (a_{bw}^m-1+\alpha+\epsilon) \times (a_{bw}^m+1-\alpha+\epsilon)$, $\epsilon$,$\alpha\in\mathbb{R}$. Let $f_{\alpha_0}(\epsilon)=f(\epsilon,\alpha_0)$ and $f_{\epsilon_0}(\alpha)=f(\epsilon_0,\alpha)$ for $\alpha_0,\epsilon_0\in\mathbb{R}$. If $\alpha_0\in[0,1]$, then $f_{\alpha_0}$ is a decreasing function for $\epsilon\in[0,a_{bw}^m-\sqrt{2}]$ and if $\epsilon_0\in[0,a_{bw}^m-\sqrt{2}]$, then $f_{\epsilon_0}$ is an increasing function for $\alpha\in[0,1]$.
	\end{lemma}
	\begin{proof}
		Let $\alpha_0\in[0,1]$. Since $0<(a_{bw}^m-1+\alpha_0-\epsilon)$ and $0<(a_{bw}^m+1-\alpha_0-\epsilon)$ for $\epsilon_0\in[0,a_{bw}^m-\sqrt{2}]$, $f_{\alpha_0}$ is a (strictly) decreasing function in this domain.\\\\
		Let $\epsilon_0\in[0,a_{bw}^m-\sqrt{2}]$. Then $f'_{\epsilon_0}(\alpha)=2(a_{bw}^m-1+\alpha-\epsilon_0)(a_{bw}^m+1-\alpha-\epsilon_0)(2-2\alpha)-(2-2\alpha).$ So, $f'_{\epsilon_0}(\alpha)=0$ gives $\alpha=1,1\pm\sqrt{(a_{bw}^m-\epsilon_0)^2-\frac{1}{2}}$. Since $f'_{\epsilon_0}(\alpha)\geq0$ for $x\in[1-\sqrt{(a_{bw}^m-\epsilon_0)^2-\frac{1}{2}},1]$, $f_{\epsilon_0}$ is an increasing function in this interval. As $0\leq\epsilon_0\leq a_{bw}^m-\sqrt{2}$, we get 
		$(a_{bw}^m-\epsilon_0)^2-\frac{1}{2}\geq(2-a_{bw}^m+\sqrt{2})^2-\frac{1}{2}=\frac{3}{2}$. This implies $1-\sqrt{(a_{bw}^m-\epsilon_0)^2-\frac{1}{2}}<0.$ So, $f_{\epsilon_0}$ is an increasing function in $[0,1]$.
	\end{proof}
	
	\begin{proposition}	\label{max_cor_value_4}
		Let $i\in\{1,2,...,n\}\setminus\{b,w\}$ be such that $\tilde{a}_{bi}=\tilde{a}_{iw}=\tilde{a}_{bw}$, and let $\epsilon_\alpha$ be the CV of $
		a_{bw}^l(\alpha) \times a_{bw}^u(\alpha) \times a_{bw}^u(\alpha) \times a_{bw}^l(\alpha) > a_{bw}^l(\alpha) \times a_{bw}^u(\alpha)$. Then $\displaystyle\Sup_{\alpha\in[0,1]}\{\epsilon_\alpha\} = \epsilon_1$. Consequently, $\epsilon_1$ is the maximum possible CV for Subcase 1 of Case 3. Also $\epsilon_1\leq$ CI. 
	\end{proposition}
	\begin{proof}
		If $\tilde{a}_{bw}=\tilde{9}$, then $\epsilon_\alpha=5.2279$ for all $\alpha\in[0,1]$.\\\\
		Now assume that $\tilde{a}_{bw}\neq\tilde{9}$. Let $\alpha\in[0,1)$. Suppose, if possible, $\epsilon_1<\epsilon_\alpha$. By Lemma \ref{upper_bound} , $\epsilon_\alpha,\epsilon_1<a_{bw}^m-\sqrt{2}$. As $\epsilon_\alpha$ is the CV of given inequality, $f_\alpha(\epsilon_\alpha)=0$. So, by Lemma \ref{dif_fun} , we get $0=f_\alpha(\epsilon_\alpha)<f_\alpha(\epsilon_1)\leq f_1(\epsilon_1)=0$ which is contradiction. Hence $\epsilon_\alpha\leq\epsilon_1$. So, in both cases, we get $\displaystyle\Sup_{\alpha\in[0,1]}\{\epsilon_\alpha\} = \epsilon_1$. Now Proposition \ref{consistency_lower_upper_lower_upper_inequality} and Proposition \ref{lower_bound_3_1} implies that $\epsilon_1$ is the maximum possible CV in Case 3 and $\epsilon_1\leq$ CI respectively. Hence the proof.
	\end{proof}
	The values of such $\epsilon_1$ corresponding to different values of $\tilde{a}_{bw}$ are given in Table \ref{consistency_indices} .\\\\
	\textbf{Subcase 2.} $a_{bi_1}^l(\alpha) \times a_{i_1w}^u(\alpha) \times a_{bi_2}^u(\alpha) \times a_{i_2w}^l(\alpha) < a_{bw}^l(\alpha) \times a_{bw}^u(\alpha)$. Here the CV of $a_{bi_1}^l(\alpha) \times a_{i_1w}^u(\alpha) \times a_{bi_2}^u(\alpha) \times a_{i_2w}^l(\alpha) < a_{bw}^l(\alpha) \times a_{bw}^u(\alpha)$, say $\epsilon$, is defined as 
	\begin{multline*}
		\epsilon=\inf\{\max\{\epsilon_1,\epsilon_2,\epsilon_3,\epsilon_4,\epsilon_5,\epsilon_6\}: (a_{bi_1}^l(\alpha)\pm \epsilon_1 )\times (a_{i_1w}^u(\alpha)\pm\epsilon_2)\times(a_{bi_2}^u(\alpha)\pm\epsilon_3) \times (a_{i_2w}^l(\alpha)\pm\epsilon_4)=\\(a_{bw}^l(\alpha)\pm\epsilon_5) \times (a_{bw}^u(\alpha)\pm\epsilon_6) \text{ and }\epsilon_1,\epsilon_2,\epsilon_3,\epsilon_4,\epsilon_5,\epsilon_6\geq 0\}.
	\end{multline*}
	Consider the equation
	\begin{eqnarray}
		\label{cor_equation_1}
		(a_{bi_1}^l(\alpha)+x) \times (a_{i_1w}^u(\alpha)+x)
		\times (a_{bi_2}^u(\alpha) +x) \times (a_{i_2w}^l(\alpha)+x) =(a_{bw}^l(\alpha)-x) \times (a_{bw}^u(\alpha)-x).
	\end{eqnarray}
	It can be proven that the above equation has a positive root, $\epsilon$ is the smallest positive root of equation (\ref{cor_equation_1}) and $\epsilon$ is strictly less than $a_{bw}^l(\alpha)$ as well as $a_{bw}^u(\alpha)$.
	\begin{proposition}\label{lower_bound_3_2}
		Let $i_1,i_2\in\{1,2,...,n\}\setminus\{b,w\}$ and $\alpha\in [0,1]$ such that
		$a_{bi_1}^l(\alpha) \times a_{i_1w}^u(\alpha) \times a_{bi_2}^u(\alpha) \times a_{i_2w}^l(\alpha) < a_{bw}^l(\alpha) \times a_{bw}^u(\alpha)$, let $\epsilon$ be the CV of this inequality, and let $\epsilon^*$ be the optimal objective value in (\ref{general_minimization}). Then $\epsilon\leq\epsilon^*$. 
	\end{proposition}
	\begin{proposition}
		Let $i\in\{1,2,...,n\}\setminus\{b,w\}$ be such that $\tilde{a}_{bi}=\tilde{a}_{iw}=\tilde{1}$, and let $\epsilon_\alpha$ be the CV of $1^l(\alpha) \times 1^u(\alpha) \times 1^u(\alpha) \times 1^l(\alpha) < a_{bw}^l(\alpha) \times a_{bw}^u(\alpha)$. Then $\displaystyle\Sup_{\alpha\in[0,1]}\{\epsilon_\alpha\} = \epsilon_1$. Also, $\epsilon_1$ is the maximum possible CV for Subcase 2 of Case 3 and $\epsilon_1\leq$ CI.
	\end{proposition}
	The values of such $\epsilon_1$ corresponding to different values of $\tilde{a}_{bw}$ are given in Table \ref{consistency_indices} .\\\\
	\textbf{Case 4.} Suppose that (\ref{consistency_lower_upper_inequality}) does not hold. So, there exist $i\in\{1,2,...,n\}\setminus\{b,w\}$ and $\alpha\in[0,1]$ such that $a_{bi}^l(\alpha) \times a_{iw}^u(\alpha) > a_{bw}^u(\alpha)$. Then the CV of $a_{bi}^l(\alpha) \times a_{iw}^u(\alpha) > a_{bw}^u(\alpha)$, say $\epsilon$, is defined as 
	\begin{multline*}
		\epsilon=\inf\{\max\{\epsilon_1,\epsilon_2,\epsilon_3\}: (a_{bi_1}^l(\alpha)\pm \epsilon_1 )\times (a_{i_1w}^u(\alpha)\pm\epsilon_2) \leq (a_{bw}^u(\alpha)\pm\epsilon_3) \text{ and }\epsilon_1,\epsilon_2,\epsilon_3\geq 0\}.
	\end{multline*}
	Consider the quadratic equation
	\begin{eqnarray}
		(a_{bi}^l(\alpha) - x) \times (a_{iw}^u(\alpha) - x) = (a_{bw}^u(\alpha) + x). \nonumber
	\end{eqnarray} 
	It can be proved that there exist $0<c$ such that $c$ is a root of the above equation, $c<a_{bi}^l(\alpha)$ and $c<a_{iw}^u(\alpha)$. Note that $\epsilon$ is the smallest positive root of the above equation. So, it follows that $\epsilon$ is strictly less than $a_{bi}^l(\alpha)$ as well as $a_{iw}^u(\alpha)$.\\\\
	The above equation can be rewritten as
	\begin{eqnarray*}
		x^2- (a_{bi}^l(\alpha)+a_{iw}^u(\alpha)+1)\times x +a_{bi}^l(\alpha) \times a_{iw}^u(\alpha)-a_{bw}^u(\alpha) =0. 
	\end{eqnarray*}
	Roots of this quadratic equation are $\frac{(a_{bi}^l(\alpha)+a_{iw}^u(\alpha)+1)\pm\sqrt{\Delta}}{2},$	where
	$\Delta=(a_{bi}^l(\alpha)+a_{iw}^u(\alpha)+1)^2-4\times (a_{bi}^l(\alpha) \times a_{iw}^u(\alpha)-a_{bw}^u(\alpha)).$\\\\
	Since the above equation has a real root which is less than $a_{bi}^l(\alpha)$ and $a_{iw}^u(\alpha)$, we have $\Delta> 0$. We also have $a_{bi}^l(\alpha) \times a_{iw}^u(\alpha) > a_{bw}^u(\alpha)$. So, $a_{bi}^l(\alpha) \times a_{iw}^u(\alpha) - a_{bw}^u(\alpha)>0$. This implies $0<(a_{bi}^l(\alpha)+a_{iw}^u(\alpha)+1)^2-4(a_{bi}^l(\alpha) \times a_{iw}^u(\alpha)-a_{bw}^u(\alpha))< (a_{bi}^l(\alpha)+a_{iw}^u(\alpha)+1)^2$. So, $\sqrt{\Delta}<a_{bi}^l(\alpha)+a_{iw}^u(\alpha)+1$. Thus, both roots of the above equation are positive. Since $\epsilon$ is the smallest positive root, we get 
	\begin{eqnarray}\label{cor_value_lower_upper_inequality}
		\epsilon=\frac{(a_{bi}^l(\alpha)+a_{iw}^u(\alpha)+1)-\sqrt{\Delta}}{2},
	\end{eqnarray}
	where $\Delta=(a_{bi}^l(\alpha)+a_{iw}^u(\alpha)+1)^2-4(a_{bi}^l(\alpha) \times a_{iw}^u(\alpha)-a_{bw}^u(\alpha))$.
	\begin{proposition}\label{lower_bound_4}
		Let $i\in\{1,2,...,n\}\setminus\{b,w\}$ and $\alpha\in [0,1]$ such that
		$a_{bi}^l(\alpha) \times a_{iw}^u(\alpha) > a_{bw}^u(\alpha)$, let $\epsilon$ be the CV of this inequality, and let $\epsilon^*$ be the optimal objective value in (\ref{general_minimization}). Then $\epsilon\leq\epsilon^*$.
	\end{proposition}
	\begin{proposition}
		\label{consistency_lower_inequality}
		Let $i_1,i_2\in\{1,2,...,n\}\setminus\{b,w\}$ and $\alpha\in[0,1]$ be such that
		\begin{eqnarray}
			a_{bi_1}^l(\alpha) \leq a_{bi_2}^l(\alpha) , \quad a_{i_1w}^u(\alpha) \leq a_{i_2w}^u(\alpha),  \nonumber
		\end{eqnarray}
		and let $\epsilon_1,\epsilon_2$ be the CVs of
		\begin{eqnarray}
			a_{bi_1}^l(\alpha) \times a_{i_1w}^u(\alpha) \geq a_{bw}^u(\alpha), \quad 
			a_{bi_2}^l(\alpha) \times a_{i_2w}^u(\alpha) \geq a_{bw}^u(\alpha) \nonumber 
		\end{eqnarray}
		respectively. Then $\epsilon_1 \leq \epsilon_2$.
	\end{proposition}
	From Proposition \ref{consistency_lower_inequality} , it follows that in this case, CV is maximum for the scale \ref{tab:table1} if $\tilde{a}_{bi}=\tilde{a}_{iw}=\tilde{a}_{bw}$ for some $i\in\{1,2,...,n\}\setminus\{b,w\}.$
	\begin{proposition}
		\label{max_cor_value_2}
		Let $i\in\{1,2,...,n\}\setminus\{b,w\}$ be such that $\tilde{a}_{bi}=\tilde{a}_{iw}=\tilde{a}_{bw}$, and let $\epsilon_\alpha$ be the CV of $a_{bi}^l(\alpha) \times a_{iw}^u(\alpha) > a_{bw}^u(\alpha)$. Then $\displaystyle\Sup_{\alpha\in[0,1]}\{\epsilon_\alpha\} = \epsilon_1$. Consequently, $\epsilon_1$ is the maximum possible CV for Case 4. Also $\epsilon_1\leq$ CI.
	\end{proposition}
	\begin{proof}
		If $\tilde{a}_{bw}=\tilde{9}$, then by equation (\ref{cor_value_lower_upper_inequality}), we get $\epsilon_\alpha=5.2279$ for all $\alpha\in[0,1]$.\\\\
		Now assume that $\tilde{a}_{bw}\neq\tilde{9}$. Let $\alpha\in[0,1]$. Observe that $a_{bw}^l(\alpha)= a_{bw}^m-1+\alpha\ \text{and} \ a_{bw}^u(\alpha)= a_{bw}^m+1-\alpha$. So, by equation (\ref{cor_value_lower_upper_inequality}), we get 
		\[\epsilon_\alpha=\frac{(2a_{bw}^m+1)-\sqrt{\Delta_\alpha}}{2},\]
		where $\Delta_\alpha=(2 a_{bw}^m+1)^2-4((a_{bw}^{2m}-a_{bw}^m)-(1-\alpha)^2-(1-\alpha))$.\\\\
		Note that $\epsilon_{\alpha_1}\leq\epsilon_{\alpha_2}$ iff $\Delta_{\alpha_1}\geq\Delta_{\alpha_2}$. Since $0\leq\alpha\leq 1$, $(1-\alpha)\geq 0$. This gives $\Delta_\alpha=(2a_{bw}^m+1)^2-4((a_{bw}^m)^2-a_{bw}^m)+4(1-\alpha)^2+4(1-\alpha)) \geq (2a_{bw}^m+1)^2-4((a_{bw}^m)^2-a_{bw}^m)=\Delta_1$. Thus, $\epsilon_\alpha\leq\epsilon_1$. So, in both cases, we have $\displaystyle\Sup_{\alpha\in[0,1]}\{\epsilon_\alpha\} = \epsilon_1$.\\
		Now Proposition \ref{consistency_lower_inequality} and Proposition \ref{lower_bound_4} implies that $\epsilon_1$ is the maximum possible CV in Case 4 and $\epsilon_1\leq$ CI respectively. Hence the proof.
	\end{proof}
	For different values of $a_{bw}$, corresponding such $\epsilon_1$ are given in Table \ref{consistency_indices} .\\\\
	\textbf{Case 5.} Suppose that (\ref{consistency_upper_lower_inequality}) does not hold. So, there exist $i\in\{1,2,...,n\}\setminus\{b,w\}$ and $\alpha\in[0,1]$ such that $a_{bi}^u(\alpha) \times a_{iw}^l(\alpha) > a_{bw}^u(\alpha)$. Here, the CV of $a_{bi}^u(\alpha) \times a_{iw}^l(\alpha) > a_{bw}^u(\alpha)$, say $\epsilon$, is defined as 
	\begin{multline*}
		\epsilon=\inf\{\max\{\epsilon_1,\epsilon_2,\epsilon_3\}: (a_{bi_1}^u(\alpha)\pm \epsilon_1 )\times (a_{i_1w}^l(\alpha)\pm\epsilon_2) \leq (a_{bw}^u(\alpha)\pm\epsilon_3) \text{ and }\epsilon_1,\epsilon_2,\epsilon_3\geq 0\}.
	\end{multline*} 
	Similar to Case 4, it can be proven that $\epsilon$ is the smallest positive root of the equation $(a_{bi}^u(\alpha)-x) \times (a_{iw}^l(\alpha)-x) = (a_{bw}^l(\alpha)+x)$. So,
	\begin{eqnarray}\label{cor_value_upper_lower_inequality}
		\epsilon=\frac{(a_{bi}^u(\alpha)+a_{iw}^l(\alpha)+1)-\sqrt{\Delta}}{2},
	\end{eqnarray}
	where $\Delta=(a_{bi}^u(\alpha)+a_{iw}^l(\alpha)+1)^2-4(a_{bi}^u(\alpha) \times a_{iw}^l(\alpha)-a_{bw}^u(\alpha))$.
	\begin{proposition}\label{lower_bound_5}
		Let $i\in\{1,2,...,n\}\setminus\{b,w\}$ and $\alpha\in [0,1]$ be such that
		$a_{bi}^u(\alpha) \times a_{iw}^l(\alpha) > a_{bw}^u(\alpha)$, let $\epsilon$ be the CV of this inequality, and let $\epsilon^*$ be the optimal objective value in (\ref{general_minimization}). Then $\epsilon\leq\epsilon^*$.
	\end{proposition}
	Similar to Case 4, it can be proven that CV of $a_{bw}^u(\alpha)\times a_{bw}^l(\alpha)>a_{bw}^u(\alpha)$ at $\alpha=1$ is the maximum possible CV in Case 5 which is same as Case 4. These values are given in Table \ref{consistency_indices} .\\\\ 
	\textbf{Case 6.} Suppose that (\ref{consistency_lower_lower_inequality}) does not hold. So, there exist  $i\in\{1,2,...,n\}\setminus\{b,w\}$ and $\alpha\in[0,1]$ such that $a_{bi}^l(\alpha) \times a_{iw}^l(\alpha) > a_{bw}^l(\alpha)$. Here, the CV of $a_{bi}^l(\alpha) \times a_{iw}^l(\alpha) > a_{bw}^l(\alpha)$, say $\epsilon$, is defined as 
	\begin{multline*}
		\epsilon=\inf\{\epsilon_1,\epsilon_2,\epsilon_3: (a_{bi_1}^l(\alpha)\pm \epsilon_1 )\times (a_{i_1w}^l(\alpha)\pm\epsilon_2) \leq (a_{bw}^l(\alpha)\pm\epsilon_3) \text{ for some }\epsilon_1,\epsilon_2,\epsilon_3\geq 0\}.
	\end{multline*} 
	By similar argument as in Case 4, it can be proven that $\epsilon$ is the smallest postiove root of the equation $(a_{bi}^l(\alpha)-x) \times (a_{iw}^l(\alpha)-x) = (a_{bw}^l(\alpha)+x)$. So,
	\begin{eqnarray}\label{cor_value_lower_lower_inequality}
		\epsilon=\frac{(a_{bi}^l(\alpha)+a_{iw}^l(\alpha)+1)-\sqrt{\Delta}}{2},
	\end{eqnarray}
	where $\Delta=(a_{bi}^l(\alpha)+a_{iw}^l(\alpha)+1)^2-4(a_{bi}^l(\alpha) \times a_{iw}^l(\alpha)-a_{bw}^l(\alpha))$.
	\begin{proposition}\label{lower_bound_6}
		Let $i\in\{1,2,...,n\}\setminus\{b,w\}$ and $\alpha\in [0,1]$ such that
		$a_{bi}^l(\alpha) \times a_{iw}^l(\alpha) > a_{bw}^l(\alpha)$, let $\epsilon$ be the CV of this inequality, and let $\epsilon^*$ be the optimal objective value in (\ref{general_minimization}). Then $\epsilon\leq\epsilon^*$.
	\end{proposition}
	\begin{proposition}
		\label{consistency_lower_lower_lower_inequality}
		Let $i_1,i_2\in\{1,2,...,n\}\setminus\{b,w\}$ and $\alpha\in[0,1]$ be such that
		\begin{eqnarray}
			a_{bi_1}^l(\alpha) \leq a_{bi_2}^l(\alpha) , \quad a_{i_1w}^l(\alpha) \leq a_{i_2w}^l(\alpha),  \nonumber
		\end{eqnarray}
		and let $\epsilon_1,\epsilon_2$ be the CVs of
		\begin{eqnarray}
			a_{bi_1}^l(\alpha) \times a_{i_1w}^l(\alpha) \geq a_{bw}^l(\alpha), \quad 
			a_{bi_2}^l(\alpha) \times a_{i_2w}^l(\alpha) \geq a_{bw}^l(\alpha) \nonumber 
		\end{eqnarray}
		respectively. Then $\epsilon_1 \leq \epsilon_2$.
	\end{proposition}
	From Proposition \ref{consistency_lower_lower_lower_inequality} , it follows that in this case, CV is maximum for the scale \ref{tab:table1} if $\tilde{a}_{bi}=\tilde{a}_{iw}=\tilde{a}_{bw}$ for some $i\in\{1,2,...,n\}\setminus\{b,w\}.$
	\begin{proposition}
		\label{max_cor_value_3}
		Let $i\in\{1,2,...,n\}\setminus\{b,w\}$ be such that $\tilde{a}_{bi}=\tilde{a}_{iw}=\tilde{a}_{bw}$, and let $\epsilon_\alpha$ be the CV of $a_{bi}^l(\alpha) \times a_{iw}^l(\alpha) > a_{bw}^l(\alpha)$. Then $\displaystyle\Sup_{\alpha\in[0,1]}\{\epsilon_\alpha\} = \epsilon_1$. Consequently, $\epsilon_1$ is the maximum possible CV for Case 6. Also $\epsilon_1\leq$ CI.
	\end{proposition}
	\begin{proof}
		By equation (\ref{cor_value_lower_lower_inequality}), we have		\[\epsilon_\alpha=\frac{(a_{bw}^l(\alpha)+a_{bw}^l(\alpha)+1)-\sqrt{\Delta_\alpha}}{2},\]
		where $\Delta_\alpha=(a_{bw}^l(\alpha)+a_{bw}^l(\alpha)+1)^2-4(a_{bw}^l(\alpha) \times a_{bw}^l(\alpha)-a_{bw}^l(\alpha)) = 8a_{bw}^l(\alpha)+1$. Since $\tilde{a}_{bw}\in\{\tilde{2},\tilde{3},...,\tilde{9}\},$ we get $1\leq a_{bw}^l(\alpha)\leq 9$ for all $\alpha\in[0,1]$. Define $f(x)=(2x+1)-\sqrt{8x+1},x\geq 0$. We know that if $\alpha_1\leq \alpha_2$, then $a_{bw}^l(\alpha_1)\leq a_{bw}^l(\alpha_2)$. So, to prove $\epsilon_\alpha\leq\epsilon_1$ for all $\alpha\in[0,1]$, it is sufficient to prove that $f$ is an increasing function in $[1,9]$. Here $f'(x)=2-\frac{4}{\sqrt{8x+1}}$. So, $f'(x)=0$ implies $x=\frac{3}{8}$. Observe that $f'(x)<0 \ \text{for}\  x\in[0,\frac{3}{8})$ and $f'(x)>0 \ \text{for}\ x\in(\frac{3}{8},\infty)$. Thus, $f$ is an increasing function in the domain $(\frac{3}{8},\infty)$, particularly in $[1,9]$. This gives $\displaystyle\Sup_{\alpha\in[0,1]}\{\epsilon_\alpha\} = \epsilon_1$. \\
		Now Proposition \ref{consistency_lower_lower_lower_inequality} and Proposition \ref{lower_bound_6} implies that $\epsilon_1$ is the maximum possible CV in Case 6 and $\epsilon_1\leq$ CI respectively. Hence the proof.
	\end{proof}
	For different values of $a_{bw}$, corresponding such $\epsilon_1$ are given in Table \ref{consistency_indices} .\\\\
	\textbf{Case 7.} Suppose that (\ref{consistency_increasing_1}) does not hold. So, there exist $i \in\{1,2,...,n\}\setminus\{b,w\}$ and $\alpha_1,\alpha_2\in[0,1]$ such that $\alpha_1<\alpha_2$ and $\frac{a_{bi}^l(\alpha_2) \times a_{iw}^u(\alpha_2)}{a_{bw}^u(\alpha_2)}<\frac{a_{bi}^l(\alpha_1) \times a_{iw}^u(\alpha_1)}{a_{bw}^u(\alpha_1)}$. Here the CV of this inequality, say $\epsilon$, is defined as 
	\begin{multline*}
		\epsilon=\inf\{\max\{\epsilon_1,\epsilon_2,\epsilon_3,\epsilon_4,\epsilon_5,\epsilon_6\}: \frac{(a_{bi}^l(\alpha_1)\pm \epsilon_1) \times (a_{iw}^u(\alpha_1)\pm \epsilon_2)}{(a_{bw}^u(\alpha_1)\pm \epsilon_3)}\leq\frac{(a_{bi}^l(\alpha_2)\pm \epsilon_4) \times (a_{iw}^u(\alpha_2)\pm \epsilon_5)}{(a_{bw}^u(\alpha_2)\pm \epsilon_6)}\\ \text{ and }\epsilon_1,\epsilon_2,\epsilon_3,\epsilon_4,\epsilon_5,\epsilon_6\geq 0\}.
	\end{multline*} 
	So, $\epsilon$ is the smallest positive root of the equation 
	\begin{eqnarray*}
		(a_{bi}^l(\alpha_2) +x)\times (a_{iw}^u(\alpha_2)+x) \times (a_{bw}^u(\alpha_1)+x) = (a_{bi}^l(\alpha_1)-x) \times (a_{iw}^u(\alpha_1)-) \times (a_{bw}^u(\alpha_2)-x).
	\end{eqnarray*}	
	\begin{proposition}\label{lower_bound_7}
		Let $i\in\{1,2,...,n\}\setminus\{b,w\}$, $\alpha_1,\alpha_2\in [0,1]$ such that $\alpha_1<\alpha_2$ and	$\frac{a_{bi}^l(\alpha_2) \times a_{iw}^u(\alpha_2)}{a_{bw}^u(\alpha_2)}<\frac{a_{bi}^l(\alpha_1) \times a_{iw}^u(\alpha_1)}{a_{bw}^u(\alpha_1)}$, let $\epsilon$ be the CV of this inequality, and let $\epsilon^*$ be the optimal objective value in (\ref{general_minimization}). Then $\epsilon\leq\epsilon^*$.
	\end{proposition}
	\begin{proposition}\label{upper_bound_1}
		Let $i \in\{1,2,...,n\}\setminus\{b,w\}$ and $\alpha_1,\alpha_2\in[0,1]$ be such that $\alpha_1<\alpha_2$ and $\frac{a_{bi}^l(\alpha_2) \times a_{iw}^u(\alpha_2)}{a_{bw}^u(\alpha_2)}<\frac{a_{bi}^l(\alpha_1) \times a_{iw}^u(\alpha_1)}{a_{bw}^u(\alpha_1)}$, and let $\epsilon$ be the CV of this inequality. Then $$\epsilon\leq \max\biggl\{\frac{a_{bi}^l(\alpha_2)-a_{bi}^l(\alpha_1)}{2},\frac{a_{iw}^u(\alpha_1)-a_{iw}^u(\alpha_2)}{2},\frac{a_{bw}^u(\alpha_1)-a_{bw}^u(\alpha_2)}{2}\biggr\}.$$ 
	\end{proposition}
	\begin{proof}
		First observe that $\frac{(a_{bi}^l(\alpha_2)-a_{bi}^l(\alpha_1))}{2}$, $\frac{(a_{iw}^u(\alpha_1)-a_{iw}^u(\alpha_2))}{2}$, $\frac{(a_{bw}^u(\alpha_1)-a_{bw}^u(\alpha_2))}{2}\geq 0$. Also note that 
		\begin{eqnarray*}
			&&\frac{\big(a_{bi}^l(\alpha_2)-\frac{a_{bi}^l(\alpha_2)-a_{bi}^l(\alpha_1)}{2}\big) \times \big(a_{iw}^u(\alpha_2)+\frac{a_{iw}^u(\alpha_1)-a_{iw}^u(\alpha_2)}{2}\big)}{a_{bw}^u(\alpha_2)+\frac{a_{bw}^u(\alpha_1)-a_{bw}^u(\alpha_2)}{2}}\\
			&=&\frac{\frac{a_{bi}^l(\alpha_1)+a_{bi}^l(\alpha_2)}{2}\times \frac{a_{bw}^u(\alpha_1)+a_{bw}^u(\alpha_2)}{2}}{\frac{a_{bw}^u(\alpha_1)+a_{bw}^u(\alpha_2)}{2}}\\
			&=&\frac{\big(a_{bi}^l(\alpha_1)+\frac{a_{bi}^l(\alpha_2)-a_{bi}^l(\alpha_1)}{2}\big) \times \big(a_{iw}^u(\alpha_1)-\frac{a_{iw}^u(\alpha_1)-a_{iw}^u(\alpha_2)}{2}\big)}{a_{bw}^u(\alpha_1)-\frac{a_{bw}^u(\alpha_1)-a_{bw}^u(\alpha_2)}{2}}.
		\end{eqnarray*}
		So, by definition of CV, we get $\epsilon\leq\max\biggl\{\frac{a_{bi}^l(\alpha_2)-a_{bi}^l(\alpha_1)}{2},\frac{a_{iw}^u(\alpha_1)-a_{iw}^u(\alpha_2)}{2},\frac{a_{bw}^u(\alpha_1)-a_{bw}^u(\alpha_2)}{2}\biggr\}.$
	\end{proof}
	\begin{remark}
		For the scale \ref{tab:table1} , we have $a_{bi}^l(\alpha)=\begin{cases}
			a_{bi}^m, \quad \quad \quad\quad\ \ \text{if }\tilde{a}_{bi}=\tilde{1} \text{ or } \tilde{9}\\ 
			a_{bi}^m-1+\alpha, \quad \text{otherwise} 
		\end{cases}$ and $a_{iw}^l(\alpha)=\begin{cases}
			a_{iw}^m, \quad \quad \quad\quad\ \ \text{if }\tilde{a}_{iw}=\tilde{1} \text{ or } \tilde{9}\\ 
			a_{iw}^m+1-\alpha, \quad \text{otherwise}
		\end{cases}$ for all $i\in\{1,2,...,n\}$ and $\alpha\in[0,1]$. So, we get 
		\begin{eqnarray}\label{upper_value}
			\epsilon\leq\max\biggl\{\frac{a_{bi}^l(\alpha_2)-a_{bi}^l(\alpha_1)}{2},\frac{a_{iw}^u(\alpha_1)-a_{iw}^u(\alpha_2)}{2},\frac{a_{bw}^u(\alpha_1)-a_{bw}^u(\alpha_2)}{2}\biggr\}\leq \frac{\alpha_2-\alpha_1}{2}\leq 0.5.
		\end{eqnarray}
	\end{remark}
	
	\textbf{Case 8.} Suppose that (\ref{consistency_increasing_2}) does not hold. So, there exist $i \in\{1,2,...,n\}\setminus\{b,w\}$ and $\alpha_1,\alpha_2\in[0,1]$ such that $\alpha_1<\alpha_2$ and $\frac{a_{bi}^u(\alpha_2) \times a_{iw}^l(\alpha_2)}{a_{bw}^u(\alpha_2)}<\frac{a_{bi}^u(\alpha_1) \times a_{iw}^l(\alpha_1)}{a_{bw}^u(\alpha_1)}$. Here, the CV of this inequality, say $\epsilon$, is defined as 
	\begin{multline*}
		\epsilon=\inf\{\max\{\epsilon_1,\epsilon_2,\epsilon_3,\epsilon_4,\epsilon_5,\epsilon_6\}: \frac{(a_{bi}^u(\alpha_1)\pm \epsilon_1) \times (a_{iw}^l(\alpha_1)\pm \epsilon_2)}{(a_{bw}^u(\alpha_1)\pm \epsilon_3)}\leq\frac{(a_{bi}^u(\alpha_2)\pm \epsilon_4) \times (a_{iw}^l(\alpha_2)\pm \epsilon_5)}{(a_{bw}^u(\alpha_2)\pm \epsilon_6)}\\ \text{ and }\epsilon_1,\epsilon_2,\epsilon_3,\epsilon_4,\epsilon_5,\epsilon_6\geq 0\}.
	\end{multline*}
	So, $\epsilon$ is the smallest positive root of the equation
	\begin{eqnarray*}
		(a_{bi}^u(\alpha_2) +x)\times (a_{iw}^l(\alpha_2)+x) \times (a_{bw}^u(\alpha_1)+x)=(a_{bi}^u(\alpha_1)-x) \times (a_{iw}^l(\alpha_1)-x) \times (a_{bw}^u(\alpha_2)-x).
	\end{eqnarray*}
	\begin{proposition}\label{lower_bound_8}
		Let $i\in\{1,2,...,n\}\setminus\{b,w\}$, $\alpha_1,\alpha_2\in [0,1]$ such that $\alpha_1<\alpha_2$ and	$\frac{a_{bi}^u(\alpha_2) \times a_{iw}^l(\alpha_2)}{a_{bw}^u(\alpha_2)}<\frac{a_{bi}^u(\alpha_1) \times a_{iw}^l(\alpha_1)}{a_{bw}^u(\alpha_1)}$, let $\epsilon$ be the CV of this inequality, and let $\epsilon^*$ be the optimal objective value in (\ref{general_minimization}). Then $\epsilon\leq\epsilon^*$.
	\end{proposition}
	Similar to Case 7, it can be proven that $	\epsilon\leq\max\biggl\{\frac{a_{bi}^u(\alpha_1)-a_{bi}^u(\alpha_2)}{2},\frac{a_{iw}^l(\alpha_2)-a_{iw}^l(\alpha_1)}{2},\frac{a_{bw}^u(\alpha_1)-a_{bw}^u(\alpha_2)}{2}\biggr\}$, and consequently for scale \ref{tab:table1} , $\epsilon\leq \frac{\alpha_2-\alpha_1}{2}\leq 0.5.$\\\\
	\textbf{Case 9.} Suppose that (\ref{consistency_increasing_3}) does not hold. So, there exist $i \in\{1,2,...,n\}\setminus\{b,w\}$ and $\alpha_1,\alpha_2\in[0,1]$ such that $\alpha_1<\alpha_2$ and $\frac{a_{bi}^l(\alpha_2) \times a_{iw}^l(\alpha_2)}{a_{bw}^l(\alpha_2)}<\frac{a_{bi}^l(\alpha_1) \times a_{iw}^l(\alpha_1)}{a_{bw}^l(\alpha_1)}$. Here the CV of this inequality, say $\epsilon$, is defined as 
	\begin{multline*}
		\epsilon=\inf\{\max\{\epsilon_1,\epsilon_2,\epsilon_3,\epsilon_4,\epsilon_5,\epsilon_6\}: \frac{(a_{bi}^l(\alpha_1)\pm \epsilon_1) \times (a_{iw}^l(\alpha_1)\pm \epsilon_2)}{(a_{bw}^l(\alpha_1)\pm \epsilon_3)}\leq\frac{(a_{bi}^l(\alpha_2)\pm \epsilon_4) \times (a_{iw}^l(\alpha_2)\pm \epsilon_5)}{(a_{bw}^l(\alpha_2)\pm \epsilon_6)}\\ \text{ and }\epsilon_1,\epsilon_2,\epsilon_3,\epsilon_4,\epsilon_5,\epsilon_6\geq 0\}.
	\end{multline*}
	So, $\epsilon$ is the smallest positive root of the equation
	\begin{eqnarray*}
		(a_{bi}^l(\alpha_2) +x)\times (a_{iw}^l(\alpha_2)+x) \times (a_{bw}^l(\alpha_1)+x)=(a_{bi}^l(\alpha_1)-x) \times (a_{iw}^l(\alpha_1)-x) \times (a_{bw}^l(\alpha_2)-x).
	\end{eqnarray*}
	\begin{proposition}\label{lower_bound_9}
		Let $i\in\{1,2,...,n\}\setminus\{b,w\}$, $\alpha_1,\alpha_2\in [0,1]$ such that $\alpha_1<\alpha_2$ and	$\frac{a_{bi}^l(\alpha_2) \times a_{iw}^l(\alpha_2)}{a_{bw}^l(\alpha_2)}<\frac{a_{bi}^l(\alpha_1) \times a_{iw}^l(\alpha_1)}{a_{bw}^l(\alpha_1)}$, let $\epsilon$ be the CV of this inequality, and let $\epsilon^*$ be the optimal objective value in (\ref{general_minimization}). Then $\epsilon\leq\epsilon^*$.
	\end{proposition}
	Similar to Case 7, it can be proven that $	\epsilon\leq\max\biggl\{\frac{a_{bi}^l(\alpha_2)-a_{bi}^l(\alpha_1)}{2},\frac{a_{iw}^l(\alpha_2)-a_{iw}^l(\alpha_1)}{2},\frac{a_{bw}^l(\alpha_2)-a_{bw}^l(\alpha_1)}{2}\biggr\}$ and consequently for scale \ref{tab:table1}, $\epsilon\leq \frac{\alpha_2-\alpha_1}{2}\leq 0.5.$\\
	\begin{remark}
		From the above discussion, it follows that $\max\{\epsilon_i:i=1,2,...,n\}\leq$ CI, where $\epsilon_i$ is the maximum possible CV in Case i. From numerical examples, we have observed that $\max\{\epsilon_i:i=1,2,...,n\}=$ CI, but this observation lacks theoratical support. Due to this limitation, we cannot calculate the exact value of CR but we can calculate its upper bound, which is sufficient to determine whether a weight set is acceptable or not in most cases. 
	\end{remark}
	\begin{table}[ht]
		\caption{Lower bounds of Consistency Index \label{consistency_indices}}
		\centering
		\begin{tabular}{|c|c|c|c|c|c|}
			\hline
			\multirow{4}{*}{$\tilde{a}_{bw}$} & \multicolumn{4}{c|}{Maximum possible CV} & Lower \\
			\cline{2-5}
			& \multirow{3}{*}{Case 1 \& 2} & Case 3:& \multirow{3}{*}{Case 3:Subcase 2} &\multirow{3}{*}{Case 7,8 \& 9}& bounds \\
			& & Subcase 1, & \multirow{2}{*}{} && of CI\\
			&&Case 4,5 \& 6 &&&$\epsilon$\\
			\hline
			$\tilde{2}$& $0.5$ & $0.4384$ & $0.3027$ &\multirow{8}{*}{$\leq0.5$}& $0.5$\\
			\cline{1-4}
			\cline{6-6}
			$\tilde{3}$& $1$ & $1$ & $0.5615$ && $1$\\
			\cline{1-4}
			\cline{6-6}
			$\tilde{4}$& $1.5$ & $1.6277$ & $0.7912$ && $1.6277$\\
			\cline{1-4}
			\cline{6-6}
			$\tilde{5}$& $2$ & $2.2984$ & $1$ && $2.2984$\\
			\cline{1-4}
			\cline{6-6}
			$\tilde{6}$ & $2.5$ & $3$ & $1.1925$& & $3$\\
			\cline{1-4}
			\cline{6-6}
			$\tilde{7}$& $3$ & $3.725$ & $1.2722$&& $3.725$\\
			\cline{1-4}
			\cline{6-6}
			$\tilde{8}$& $3.5$ & $4.4688$ & $1.5413$ && $4.4688$\\
			\cline{1-4}
			\cline{6-6}
			$\tilde{9}$& $4$ & $5.2279$ & $1.7015$ & &$5.2279$\\
			\hline
		\end{tabular}
	\end{table}
	\subsection{Numerical examples}
	Now we discuss some numerical examples to illustrate the proposed model and compare it with FBWM\cite{guo2017fuzzy}. For that, first we introduce a notation.\\\\
	For $m\geq2$, let $F_m=\{0,\frac{1}{m-1},\frac{2}{m-1},...,\frac{m-1}{m-1}=1\}$, i.e., $F_m$ is subset of $[0,1]$ having $m$ elements that partitions $[0,1]$ into $m-1$ uniform sub-intervals. For example, $F_2=\{0,1\}$, $F_{10}=\{0,0.1,0.2,...,0.9,1\}$. Let $(\tilde{A}_b,\tilde{A}_w)$ be a triangular FPCS, let $\epsilon^*$ be the optimal objective value of problem (\ref{general_minimization}), and let $\eta_{F_m}^*$ be the optimal objective value of problem (\ref{optimization}) formulated using $F_m$. Note that $||F_m||_{\infty}$=$\frac{1}{m-1}$. So, by Theorem \ref{doa}, we get $|\epsilon^*-\eta_{F_m}^*|\leq\frac{1}{m-1}$. This implies that the DoA of an approximate weight set obtained using $F_m$ in problem (\ref{optimization}) is $\frac{1}{m-1}$. Therefore, this relationship helps to find the DoA of the weights set. We can also set the value of $m$ to obtain a weight set with the desired DoA.\\\\
	Observe that for $m=2$, problem (\ref{optimization}) is precisely the minimization problem of FBWM\cite{guo2017fuzzy}, i.e., $F_2^*$ is the optimal accuracy of FBWM. From the above discussion, it follows that an optimal weight set of FBWM has DoA equal to $1$. By taking higher value of $m$, we can get better DoA. For example, $m=11$ gives DoA equal to $0.1$. On the other hand, to get DoA equal to $0.01$, we need to take $m=101$. In the following examples, we have computed the weights by taking $m=2^i+1, i\in \mathbb{N}\cup \{0\}$, because for such $m$, the successive partition divides each sub-interval of the previous partition into two equal intervals.\\\\
	\textbf{Example 1.} Let $C=\{c_1,c_2,...,c_5\}$ be the set of decision criteria with $c_2$ and $c_5$ as the best and the worst criterion respectively. Let $\tilde{A}_b=(\tilde{2},\tilde{1},\tilde{4},\tilde{2},\tilde{8})$ and $\tilde{A}_w=(\tilde{3},\tilde{8},\tilde{5},\tilde{4},\tilde{1})^T$ be the best-to-other and the other-to-worst vectors respectively. Then computed weights for $i=0,4$ and $6$ are given in Table \ref{example_1}. \\\\	
	\textbf{Example 2.} Let $C=\{c_1,c_2,...,c_5\}$ be the set of decision criteria with $c_2$ and $c_5$ as the best and the worst criterion respectively. Let $\tilde{A}_b=(\tilde{3},\tilde{1},\tilde{3},\tilde{2},\tilde{6})$ and $\tilde{A}_w=(\tilde{2},\tilde{6},\tilde{6},\tilde{3},\tilde{1})^T$ be the best-to-other and the other-to-worst vectors respectively. Then computed weights for $i=0,4$ and $6$ are given in Table \ref{example_2}.\\\\
	\begin{table}[ht]
		\caption{Computed weights: Example 1\label{example_1}}
		\centering
		\begin{tabular}{|c|c|c|c|c|c|c|}
			\hline
			\multirow{2}{*}{Criterion}&\multicolumn{2}{c|}{$m=2$ (Same as FBWM\cite{guo2017fuzzy})}&\multicolumn{2}{c|}{$m=17$}&\multicolumn{2}{c|}{$m=129$}\\
			\cline{2-7}
			&Interval-weight&Average&Interval-weight&Average&Interval-weight&Average\\
			\hline
			$c_1$&$[0.1245,0.2106]$&$0.1676$&$[0.1269,0.2072]$&$0.1671$&$[0.1270,0.2071]$&$0.1671$\\
			$c_2$&$[0.3863,0.4780]$&$0.4322$&$[0.3943,0.4778]$&$0.4365$&$[0.3946,0.4778]$&$0.4367$\\
			$c_3$&$[0.1513,0.1907]$&$0.1710$&$[0.1516,0.1837]$&$0.1677$&$[0.1516,0.1836]$&$0.1676$\\
			$c_4$&$[0.1319,0.2477]$&$0.1898$&$[0.1336,0.2431]$&$0.1884$&$[0.1336,0.2430]$&$0.1884$\\
			$c_5$&$[0.0418,0.0522]$&$0.0470$&$[0.0420,0.0509]$&$0.0465$&$[0.0420,0.0509]$&$0.0465$\\
			\hline
			$\epsilon_{F_m}^*$&\multicolumn{2}{c|}{$1.3945$}&\multicolumn{2}{c|}{$1.3945$}&\multicolumn{2}{c|}{$1.3945$}\\
			\hline
			CR$\leq$&\multicolumn{2}{c|}{$0.3120$}&\multicolumn{2}{c|}{$0.3120$}&\multicolumn{2}{c|}{$0.3120$}\\
			\hline
			DoA&\multicolumn{2}{c|}{$1$}&\multicolumn{2}{c|}{$0.0625$}&\multicolumn{2}{c|}{$0.0078125$}\\
			\hline
		\end{tabular}
	\end{table}	
	\begin{table}[ht]
		\caption{Computed weights: Example 2\label{example_2}}
		\centering
		\begin{tabular}{|c|c|c|c|c|c|c|}
			\hline
			\multirow{2}{*}{Criterion}&\multicolumn{2}{c|}{$m=2$ (Same as FBWM\cite{guo2017fuzzy})}&\multicolumn{2}{c|}{$m=17$}&\multicolumn{2}{c|}{$m=129$}\\
			\cline{2-7}
			&Interval-weight&Average&Interval-weight&Average&Interval-weight&Average\\
			\hline
			$c_1$&$[0.0801,0.2040]$&$0.1421$&$[0.0835,0.1953]$&$0.1394$&$[0.0836,0.1951]$&$0.1394$\\
			$c_2$&$[0.3143,0.4476]$&$0.3810$&$[0.3364,0.4476]$&$0.3920$&$[0.3368,0.4476]$&$0.3922$\\
			$c_3$&$[0.2118,0.3140]$&$0.2629$&$[0.2118,0.2931]$&$0.2525$&$[0.2118,0.2928]$&$0.2523$\\
			$c_4$&$[0.1073,0.2534]$&$0.1804$&$[0.1117,0.2438]$&$0.1778$&$[0.1118,0.2436]$&$0.1777$\\
			$c_5$&$[0.0454,0.0652]$&$0.0553$&$[0.0463,0.0617]$&$0.0540$&$[0.0463,0.0617]$&$0.0540$\\
			\hline
			$\epsilon_{F_m}^*$&\multicolumn{2}{c|}{$1.5360$}&\multicolumn{2}{c|}{$1.5360$}&\multicolumn{2}{c|}{$1.5360$}\\
			\hline
			CR$\leq$&\multicolumn{2}{c|}{$0.5120$}&\multicolumn{2}{c|}{$0.5120$}&\multicolumn{2}{c|}{$0.5120$}\\
			\hline
			DoA&\multicolumn{2}{c|}{$1$}&\multicolumn{2}{c|}{$0.0625$}&\multicolumn{2}{c|}{$0.0078125$}\\
			\hline
		\end{tabular}
	\end{table}
	Several conclusions can be drawn from these examples.
	\begin{enumerate}
		\item As discussed earlier, for a higher value of $m$, we get a better DoA. This is also evident from the nested property of the approximate interval-weights of a criterion. In other words, as the value of $m$ increases, the corresponding interval of criterion weights becomes smaller. This implies that less accurate weights are removed from the interval as we consider a higher value of $m$. This shows the superiority of the proposed model over FBWM.
		\item We have calculated an upper bound of CR instead of exact value.
	\end{enumerate}	
	\section{A real-world application}
	In this section, we discuss a real-world application of the proposed model in ranking of risk factors in supply chain 4.0.\\\\
	Supply chain is an important part of any industry. It is basically a network of individuals, organizations, technologies and resources that deal with procurement of raw materials, their transformation into final products and distribution of these products. All strategic and operational decisions related to location of facilities, production of materials, inventories to store raw materials and final products, distribution of products, etc. are covered under the Supply Chain Management (SCM) \cite{ganeshan1995introduction}.\\\\
	As we know, the present age is the age of digitization as it is widely influenced by internet and digital technologies. Almost all sectors are affected by digitization. The industrial sector is no exception to this paradigm shift and as a result, supply chains are directly affected by it\cite{ivanov2019impact}. A Digital Supply Chain (DSC) is known as supply chain 4.0. It includes technologies from Industry 4.0\cite{lasi2014industry}, Smart Manufacturing (SM)\cite{zheng2018smart}, Internet of Things (IoT)\cite{li2015internet}, etc.\\\\
	The study of risk factors in supply chains has been an important topic. Many researchers have studied various risk factors like financial risks\cite{faizal2014risk}, supply risks\cite{yang2019research}, environmental risks\cite{faizal2014risk}, demand risks\cite{kara2020data}, etc. in depth. Digitization of supply chains has significantly increased their efficiency, but it has also created some new risks, such as cyber security\cite{colicchia2019managing}. Here we rank the risk factors involved in supply chain 4.0 in terms of their criticality which will help an organization in supply chain risk management. Zekhnini et al.\cite{zekhnini2021analytic} identified elevel such risk factors which they further divided into $5$ categories. These categories and risks are given in Table \ref{risk}. The pairwise comparison values between these risk factors are adopted from \cite{zekhnini2021analytic}. These values are given in Table \ref{data}. The weights of risk factors are calculated using the same method as in the numerical examples. These weights and the ranking of risk factors using them are shown in Table \ref{risk_weights}.\\\\
	The results are similar to those of numerical examples. Although the ranking of the risk factors is the same for all values of $m$, there are significant differences in the corresponding weight sets. For $m=17$ and $m=129$, we get better approximation of optimal weights than the existing model. A higher value of $m$ yields a weight set with better DoA, which is evident from the nested property of the approximate interval-weights. We measure the accuracy of the weight set with respect to the upper bound of the CR.	
	\section{Conclusion and future plan}
	The best-worst method is one of the widely accepted MCDM methods used to calculate the weights of decision criteria\cite{rezaei2015best}. To deal with ambiguity and uncertainty of decision judgments, Guo and Zhao\cite{guo2017fuzzy} proposed a model of BWM consisting of fuzzy sets, called Fuzzy BWM (FBWM). In this paper, we propose a model of FBWM using $\alpha$-cut intervals as an improvement over the existing model. To illustrate the model, we discuss two numerical examples and a real-world application. Some salient features of this paper are as follows.
	\begin{enumerate}
		\item In FBWM, optimal weights are calculated by solving a minimization problem formulated using only $0$ and $1$-cuts. In the proposed model, this problem is modified in such a way that $\alpha$-cut intervals for all $\alpha\in [0,1]$ are included in the weight calculation. Therefore, the entire shape of fuzzy comparison values is optimized. This also helps in reducing data loss.
		\item Fuzzy arithmetic operations are defined using $\alpha$-cut intervals. Therefore, the use of $\alpha$-cut intervals naturally involves the exact fuzzy arithmetic operations in the weight calculation unlike the existing model of FBWM which uses approximate operations.
		\item It turned out that calculation of optimal weights using problem (\ref{general_minimization}) is difficult due to the involvement of infinitely many constraints. Therefore, optimal weights are approximated using finite subsets of $[0,1]$ and Degree Of Approximation (DoA) of these approximate weight sets are estimated. The resulting weight set of the existing model of FBMW is also one of the approximate weight sets of the proposed model with DoA equal to $1$. Then we develop a way to obtain a weight set with a better DoA than this weight set. We can also get a weight set with the desired DoA.
		\item Optimal weights are approximated using an optimal solution of problem (\ref{optimization}) that is formulated using a finite subset $F$ of $[0,1]$. Note that it is a non-linear minimization problem. Therefore, it can lead to multiple solutions. To deal with this problem, first we establish the fact that for a given $F$, the collection of all approximate (defuzzified) weights of a criterion is an interval. Then we calculate the GLB and LUB of this interval using the problems (\ref{optimization_glb}) and (\ref{optimization_lub}) respectively. Finally, the center value of this interval is adopted as the approximate weight of the criterion.
		\item To measure the accuracy of a weight set, we first derive the necessary conditions for FPCS to be consistent. Based on these conditions, a lower bound of CI is calculated, which leads to an upper bound of CR. In most cases, this upper bound is sufficient to check the admissibility of weight set.
	\end{enumerate}
	Although the proposed model is superior than the existing model, it has several drawbacks.
	\begin{enumerate}
		\item A major drawback of the proposed model is that the calculation of optimal weights is difficult and hence we have to work with approximate weights.
		\item It is not known whether the necessary conditions of consistent FPCS given by the Theorem \ref{consistency} are sufficient. Because of this limitation, we can only calculate the lower bound of CI (and consequently the upper bound of CR) and not the exact values.
	\end{enumerate}
	The proposed model also opens up some important future directions.
	\begin{enumerate}
		\item It would be interesting to develop a method to easily solve the minimization problem (\ref{general_minimization}) as it would lead to optimal weights.
		\item Characterization of consistent FPCS is of great importance as it will help in calculating the exact values of CI and CR so that we have better estimation of accuracy of a weight set.
		\item Using a similar technique, the best-worst method can be extended to other extensions of classical sets such as intuitionistic fuzzy sets, hesitant fuzzy sets, etc.
	\end{enumerate}
	\section*{Acknowledgements}
	The first author gratefully acknowledges the JRF (NET) from CSIR, India for financial support to carry out the research work.

	\begin{table}
		\caption{Risk factors in supply chain 4.0\label{risk}\cite{zekhnini2021analytic}}
		\centering
		\begin{tabular}{|c||c|c|}
			\hline
			Sr no.&Types of risk& risks\\
			\hline
			\hline
			\multirow{2}{*}{$1$}&\multirow{2}{*}{Supply risks $(c_1)$}&Product arrival variability $(c_{11})$\\
			&&Loss of suppliers $(c_{12})$\\
			\hline
			\hline
			\multirow{3}{*}{$2$}&\multirow{3}{*}{Industry 4.0 risks $(c_2)$}&Cyber-attack $(c_{21})$\\
			&&Information security $(c_{22})$\\
			&& Computer security $(c_{23})$\\
			\hline
			\hline
			\multirow{2}{*}{$3$}&\multirow{2}{*}{Demand risks $(c_3)$}&Fluctuation $(c_{31})$\\
			&&Demand $(c_{32})$\\
			\hline
			\hline
			\multirow{2}{*}{$4$}&\multirow{2}{*}{Operational risks $(c_4)$}&Shortage of skilled workers $(c_{41})$\\
			&&Breakdown of fractioning $(c_{42})$\\
			\hline
			\hline
			\multirow{2}{*}{$5$}&\multirow{2}{*}{Financial risks $(c_5)$}&Macroeconomic fluctuation $(c_{51})$\\
			&&Coordination in supply chain $(c_{52})$\\
			\hline		
		\end{tabular}
	\end{table} 
	\begin{table}
		\caption{Best-to-other and other-to-worst vectors\label{data}\cite{zekhnini2021analytic}}
		\centering
		\begin{subtable}[h]{0.45\textwidth}\caption{Types of risk ($c_1-c_5$)}
			\begin{tabular}{|c||c|c|c|c|c|}
				\hline
				Criteria & $c_1$& $c_2$& $c_3$& $c_4$& $c_5$(worst)\\
				\hline
				\hline
				$c_1$&-&-&-&-&$\tilde{4}$\\
				\hline
				$c_2$ (best)&$\tilde{3}$&$\tilde{1}$&$\tilde{5}$&$\tilde{4}$&$\tilde{5}$\\
				\hline
				$c_3$&-&-&-&-&$\tilde{2}$\\
				\hline
				$c_2$&-&-&-&-&$\tilde{2}$\\
				\hline
				$c_5$&-&-&-&-&$\tilde{1}$\\
				\hline
			\end{tabular}
		\end{subtable}
		\begin{subtable}[h]{0.45\textwidth}\caption{Risks of $c_1$}
			\begin{tabular}{|c||c|c|}
				\hline
				Criteria & $c_{11}$& $c_{12}$ (worst)\\
				\hline
				\hline
				$c_{11}$ (best)&$\tilde{1}$&$\tilde{2}$\\
				\hline
				$c_{12}$&-&$\tilde{1}$\\
				\hline
			\end{tabular}
		\end{subtable}
		\begin{subtable}[h]{0.45\textwidth}\caption{Risks of $c_2$}
			\begin{tabular}{|c||c|c|c|}
				\hline
				Criteria & $c_{21}$& $c_{22}$& $c_{23}$ (worst)\\
				\hline
				\hline
				$c_{21}$ (best)&$\tilde{1}$&$\tilde{2}$&$\tilde{6}$\\
				\hline
				$c_{22}$&-&-&$\tilde{3}$\\
				\hline
				$c_{23}$&-&-&$\tilde{1}$\\
				\hline
			\end{tabular}
		\end{subtable}
		\begin{subtable}[h]{0.45\textwidth}\caption{Risks of $c_3$}
			\begin{tabular}{|c||c|c|}
				\hline
				Criteria & $c_{31}$& $c_{32}$ (worst)\\
				\hline
				\hline
				$c_{31}$ (best)&$\tilde{1}$&$\tilde{2}$\\
				\hline
				$c_{32}$&-&$\tilde{1}$\\
				\hline
			\end{tabular}
		\end{subtable}
		\begin{subtable}[h]{0.45\textwidth}\caption{Risks of $c_4$}
			\begin{tabular}{|c||c|c|}
				\hline
				Criteria & $c_{41}$& $c_{42}$ (worst)\\
				\hline
				\hline
				$c_{41}$ (best)&$\tilde{1}$&$\tilde{7}$\\
				\hline
				$c_{42}$&-&$\tilde{1}$\\
				\hline
			\end{tabular}
		\end{subtable}
		\begin{subtable}[h]{0.45\textwidth}\caption{Risks of $c_5$}
			\begin{tabular}{|c||c|c|}
				\hline
				Criteria & $c_{51}$& $c_{52}$ (worst)\\
				\hline
				\hline
				$c_{51}$ (best)&$\tilde{1}$&$\tilde{5}$\\
				\hline
				$c_{52}$&-&$\tilde{1}$\\
				\hline
			\end{tabular}
		\end{subtable}
	\end{table}
	\begin{sidewaystable}[h]
		\caption {Computed weights: Risk management in supply chain 4.0\label{risk_weights}}
		\centering
		\begin{adjustbox}{width=\textwidth}
			\small
			\begin{tabular}{|c|c|c|c|c|c|c|c|c|c|c|c|c|c|c|c|c|c|c|c|c|c|c|c|}
				\hline
				\multicolumn{8}{|c|}{$m=2$ (Same as FBWM\cite{guo2017fuzzy})}&\multicolumn{8}{c|}{$m=17$ }&\multicolumn{8}{c|}{$m=129$ }\\
				\hline
				\multirow{2}{*}{Type of risk}&\multirow{2}{*}{Interval-weight}&\multirow{2}{*}{Average}&\multirow{2}{*}{Risk}&\multicolumn{2}{c|}{Local weight}&\multirow{2}{*}{Global weight}&\multirow{2}{*}{Ranking}&\multirow{2}{*}{Type of risk}&\multirow{2}{*}{Interval-weight}&\multirow{2}{*}{Average}&\multirow{2}{*}{Risk}&\multicolumn{2}{c|}{Local weight}&\multirow{2}{*}{Global weight}&\multirow{2}{*}{Ranking}&\multirow{2}{*}{Type of risk}&\multirow{2}{*}{Interval-weight}&\multirow{2}{*}{Average}&\multirow{2}{*}{Risk}&\multicolumn{2}{c|}{Local weight}&\multirow{2}{*}{Global weight}&\multirow{2}{*}{Ranking}\\
				\cline{5-6}
				\cline{13-14}
				\cline{21-22}
				&&&&Interval-weight&Average&&&&&&&Interval-weight&Average&&&&&&&Interval-weight&Average&&\\
				\hline
				\multirow{4}{*}{$c_1$}&\multirow{4}{*}{$[0.2184,0.2605]$}&\multirow{4}{*}{$0.2395$}&$c_{11}$&$[0.6429,0.6667]$&$0.6548$&$0.1568$&$2$&\multirow{4}{*}{$c_1$}&\multirow{4}{*}{$[0.2190,0.2495]$}&\multirow{4}{*}{$0.2343$}&$c_{11}$&$[0.6667,0.6667]$&$0.6667$&$0.1562$&$2$&\multirow{4}{*}{$c_1$}&\multirow{4}{*}{$[0.2198,0.2494]$}&\multirow{4}{*}{$0.2346$}&$c_{11}$&$[0.6667,0.6667]$&$0.6667$&$0.1564$&$2$\\
				&&&$c_{12}$&$[0.3333,0.3571]$&$0.3452$&$0.0827$&$5$&&&&$c_{12}$&$[0.3333,0.3333]$&$0.3333$&$0.0781$&$5$&&&&$c_{12}$&$[0.3333,0.3333]$&$0.3333$&$0.0782$&$5$\\
				\cline{4-8}
				\cline{12-16}
				\cline{20-24}
				&&&$\epsilon_{F_m}^*$&\multicolumn{2}{c|}{$0$}&-&-&&&&$\epsilon_{F_m}^*$&\multicolumn{2}{c|}{$0$}&-&-&&&&$\epsilon_{F_m}^*$&\multicolumn{2}{c|}{$0$}&-&-\\
				\cline{4-8}
				\cline{12-16}
				\cline{20-24}
				&&&CR&\multicolumn{2}{c|}{$0$}&-&-&&&&CR&\multicolumn{2}{c|}{$0$}&-&-&&&&CR$\leq$&\multicolumn{2}{c|}{$0$}&-&-\\
				\hline
				\multirow{5}{*}{$c_2$}&\multirow{5}{*}{$[0.4304,0.4851]$}&\multirow{5}{*}{$0.4578$}&$c_{21}$&$[0.5838,0.6049]$&$0.5944$&$0.2721$&$1$&\multirow{5}{*}{$c_2$}&\multirow{5}{*}{$[0.4397,0.4855]$}&\multirow{5}{*}{$0.4626$}&$c_{21}$&$[0.5838,0.6049]$&$0.5944$&$0.2749$&$1$&\multirow{5}{*}{$c_2$}&\multirow{5}{*}{$[0.4405,0.4853]$}&\multirow{5}{*}{$0.4629$}&$c_{21}$&$[0.5838,0.6049]$&$0.5944$&$0.2751$&$1$\\
				&&&$c_{22}$&$[0.2951,0.3162]$&$0.3057$&$0.1399$&$3$&&&&$c_{22}$&$[0.2951,0.3162]$&$0.3057$&$0.1414$&$3$&&&&$c_{22}$&$[0.2951,0.3162]$&$0.3057$&$0.1415$&$3$\\
				&&&$c_{23}$&$[0.0979,0.1022]$&$0.1001$&$0.0458$&$8$&&&&$c_{23}$&$[0.0979,0.1022]$&$0.1001$&$0.0463$&$8$&&&&$c_{23}$&$[0.0979,0.1022]$&$0.1001$&$0.0463$&$8$\\
				\cline{4-8}
				\cline{12-16}
				\cline{20-24}
				&&&$\epsilon_{F_m}^*$&\multicolumn{2}{c|}{$0.1624$}&-&-&&&&$\epsilon_{F_m}^*$&\multicolumn{2}{c|}{$0.1624$}&-&-&&&&$\epsilon_{F_m}^*$&\multicolumn{2}{c|}{$0.1624$}&-&-\\
				\cline{4-8}
				\cline{12-16}
				\cline{20-24}
				&&&CR$\leq$&\multicolumn{2}{c|}{$0.0541$}&-&-&&&&CR$\leq$&\multicolumn{2}{c|}{$0.0541$}&-&-&&&&CR$\leq$&\multicolumn{2}{c|}{$0.0541$}&-&-\\
				\hline
				\multirow{4}{*}{$c_3$}&\multirow{4}{*}{$[0.0830,0.1196]$}&\multirow{4}{*}{$0.1013$}&$c_{31}$&$[0.6429,0.6667]$&$0.6548$&$0.0663$&$6$&\multirow{4}{*}{$c_3$}&\multirow{4}{*}{$[0.0846,0.1163]$}&\multirow{4}{*}{$0.1005$}&$c_{31}$&$[0.6667,0.6667]$&$0.6667$&$0.0670$&$6$&\multirow{4}{*}{$c_3$}&\multirow{4}{*}{$[0.0847,0.1145]$}&\multirow{4}{*}{$0.0996$}&$c_{31}$&$[0.6667,0.6667]$&$0.6667$&$0.0664$&$6$\\
				&&&$c_{32}$&$[0.3333,0.3571]$&$0.3452$&$0.0350$&$9$&&&&$c_{32}$&$[0.3333,0.3333]$&$0.3333$&$0.0335$&$9$&&&&$c_{32}$&$[0.3333,0.3333]$&$0.3333$&$0.0332$&$9$\\
				\cline{4-8}
				\cline{12-16}
				\cline{20-24}
				&&&$\epsilon_{F_m}^*$&\multicolumn{2}{c|}{$0$}&-&-&&&&$\epsilon_{F_m}^*$&\multicolumn{2}{c|}{$0$}&-&-&&&&$\epsilon_{F_m}^*$&\multicolumn{2}{c|}{$0$}&-&-\\
				\cline{4-8}
				\cline{12-16}
				\cline{20-24}
				&&&CR&\multicolumn{2}{c|}{$0$}&-&-&&&&CR&\multicolumn{2}{c|}{$0$}&-&-&&&&CR&\multicolumn{2}{c|}{$0$}&-&-\\
				\hline
				\multirow{4}{*}{$c_4$}&\multirow{4}{*}{$[0.0977,0.1576]$}&\multirow{4}{*}{$0.1277$}&$c_{41}$&$[0.8742,0.8750]$&$0.8746$&$0.1116$&$4$&\multirow{4}{*}{$c_4$}&\multirow{4}{*}{$[0.0993,0.1557]$}&\multirow{4}{*}{$0.1275$}&$c_{41}$&$[0.8750,0.8750]$&$0.8750$&$0.1116$&$4$&\multirow{4}{*}{$c_4$}&\multirow{4}{*}{$[0.0998,0.1557]$}&\multirow{4}{*}{$0.1278$}&$c_{41}$&$[0.8750,0.8750]$&$0.8750$&$0.1118$&$4$\\
				&&&$c_{42}$&$[0.1250,0.1258]$&$0.1254$&$0.0160$&$10$&&&&$c_{42}$&$[0.1250,0.1250]$&$0.1250$&$0.0159$&$10$&&&&$c_{42}$&$[0.1250,0.1250]$&$0.1250$&$0.0160$&$10$\\
				\cline{4-8}
				\cline{12-16}
				\cline{20-24}
				&&&$\epsilon_{F_m}^*$&\multicolumn{2}{c|}{$0$}&-&-&&&&$\epsilon_{F_m}^*$&\multicolumn{2}{c|}{$0$}&-&-&&&&$\epsilon_{F_m}^*$&\multicolumn{2}{c|}{$0$}&-&-\\
				\cline{4-8}
				\cline{12-16}
				\cline{20-24}
				&&&CR&\multicolumn{2}{c|}{$0$}&-&-&&&&CR&\multicolumn{2}{c|}{$0$}&-&-&&&&CR&\multicolumn{2}{c|}{$0$}&-&-\\
				\hline
				\multirow{4}{*}{$c_5$}&\multirow{4}{*}{$[0.0720,0.0816]$}&\multirow{4}{*}{$0.0768$}&$c_{51}$&$[0.8314,0.8333]$&$0.8324$&$0.0639$&$7$&\multirow{4}{*}{$c_5$}&\multirow{4}{*}{$[0.0731,0.0808]$}&\multirow{4}{*}{$0.0770$}&$c_{51}$&$[0.8333,0.8333]$&$0.8333$&$0.0641$&$7$&\multirow{4}{*}{$c_5$}&\multirow{4}{*}{$[0.0733,0.0807]$}&\multirow{4}{*}{$0.0770$}&$c_{51}$&$[0.8333,0.8333]$&$0.8333$&$0.0642$&$7$\\
				&&&$c_{52}$&$[0.1667,0.1686]$&$0.1677$&$0.0129$&$11$&&&&$c_{52}$&$[0.1667,0.1667]$&$0.1667$&$0.0128$&$11$&&&&$c_{52}$&$[0.1667,0.1667]$&$0.1667$&$0.0128$&$11$\\
				\cline{4-8}
				\cline{12-16}
				\cline{20-24}
				&&&$\epsilon_{F_m}^*$&\multicolumn{2}{c|}{$0$}&-&-&&&&$\epsilon_{F_m}^*$&\multicolumn{2}{c|}{$0$}&-&-&&&&$\epsilon_{F_m}^*$&\multicolumn{2}{c|}{$0$}&-&-\\
				\cline{4-8}
				\cline{12-16}
				\cline{20-24}
				&&&CR&\multicolumn{2}{c|}{$0$}&-&-&&&&CR&\multicolumn{2}{c|}{$0$}&-&-&&&&CR&\multicolumn{2}{c|}{$0$}&-&-\\
				\hline
				$\epsilon_{F_m}^*$&\multicolumn{2}{c|}{$1.0140$}&\multicolumn{5}{c|}{-}&$\epsilon_{F_m}^*$&\multicolumn{2}{c|}{$1.0370$}&\multicolumn{5}{c|}{-}&$\epsilon_{F_m}^*$&\multicolumn{2}{c|}{$1.0370$}&\multicolumn{5}{c|}{-}\\
				\hline
				CR$\leq$&\multicolumn{2}{c|}{$0.4412$}&\multicolumn{5}{c|}{-}&CR$\leq$&\multicolumn{2}{c|}{$0.4512$}&\multicolumn{5}{c|}{-}&CR$\leq$&\multicolumn{2}{c|}{$0.4512$}&\multicolumn{5}{c|}{-}\\
				\hline
				DoA&\multicolumn{7}{c|}{$1$}&DoA&\multicolumn{7}{c|}{$0.0625$}&DoA&\multicolumn{7}{c|}{$0.0078125$}\\
				\hline				
			\end{tabular}	
		\end{adjustbox}
	\end{sidewaystable}

\begin{thebibliography}{10}
		
		\bibitem{ali2019hesitant}
		Asif Ali and Tabasam Rashid.
		\newblock Hesitant fuzzy best-worst multi-criteria decision-making method and
		its applications.
		\newblock {\em International Journal of Intelligent Systems}, 34(8):1953--1967,
		2019.
		
		\bibitem{brans1985note}
		Jean-Pierre Brans and Ph~Vincke.
		\newblock Note—a preference ranking organisation method: (the promethee
		method for multiple criteria decision-making).
		\newblock {\em Management science}, 31(6):647--656, 1985.
		
		\bibitem{chen2000extensions}
		Chen-Tung Chen.
		\newblock Extensions of the topsis for group decision-making under fuzzy
		environment.
		\newblock {\em Fuzzy sets and systems}, 114(1):1--9, 2000.
		
		\bibitem{colicchia2019managing}
		Claudia Colicchia, Alessandro Creazza, and David~A Menachof.
		\newblock Managing cyber and information risks in supply chains: insights from
		an exploratory analysis.
		\newblock {\em Supply Chain Management: An International Journal},
		24(2):215--240, 2019.
		
		\bibitem{dong2021fuzzy}
		Jiuying Dong, Shuping Wan, and Shyi-Ming Chen.
		\newblock Fuzzy best-worst method based on triangular fuzzy numbers for
		multi-criteria decision-making.
		\newblock {\em Information Sciences}, 547:1080--1104, 2021.
		
		\bibitem{faizal2014risk}
		K~Faizal, PLK Palaniappan, et~al.
		\newblock Risk assessment and management in supply chain.
		\newblock {\em Global Journal of Research in Engineering: G Industrail
			Engineering}, 14(2):18--30, 2014.
		
		\bibitem{ganeshan1995introduction}
		Ram Ganeshan.
		\newblock An introduction to supply chain management.
		\newblock {\em http://lcm. csa. iisc. ernet. in/scm/supply\_chain\_intro.
			html}, 1995.
		
		\bibitem{guo2017fuzzy}
		Sen Guo and Haoran Zhao.
		\newblock Fuzzy best-worst multi-criteria decision-making method and its
		applications.
		\newblock {\em Knowledge-Based Systems}, 121:23--31, 2017.
		
		\bibitem{hwang1981methods}
		Ching-Lai Hwang and Kwangsun Yoon.
		\newblock Methods for multiple attribute decision making.
		\newblock In {\em Multiple attribute decision making}, pages 58--191. Springer,
		1981.
		
		\bibitem{ivanov2019impact}
		Dmitry Ivanov, Alexandre Dolgui, and Boris Sokolov.
		\newblock The impact of digital technology and industry 4.0 on the ripple
		effect and supply chain risk analytics.
		\newblock {\em International journal of production research}, 57(3):829--846,
		2019.
		
		\bibitem{kara2020data}
		Merve~Er Kara, Seniye {\"U}mit~Oktay F{\i}rat, and Abhijeet Ghadge.
		\newblock A data mining-based framework for supply chain risk management.
		\newblock {\em Computers \& Industrial Engineering}, 139:105570, 2020.
		
		\bibitem{klir1996fuzzy}
		George~J Klir and Bo~Yuan.
		\newblock Fuzzy sets and fuzzy logic: theory and applications.
		\newblock {\em Possibility Theory versus Probab. Theory}, 32(2):207--208, 1996.
		
		\bibitem{lasi2014industry}
		Heiner Lasi, Peter Fettke, Hans-Georg Kemper, Thomas Feld, and Michael
		Hoffmann.
		\newblock Industry 4.0.
		\newblock {\em Business \& information systems engineering}, 6:239--242, 2014.
		
		\bibitem{li2015internet}
		Shancang Li, Li~Da Xu, and Shanshan Zhao.
		\newblock The internet of things: a survey.
		\newblock {\em Information systems frontiers}, 17:243--259, 2015.
		
		\bibitem{mikhailov2004evaluation}
		Ludmil Mikhailov and Petco Tsvetinov.
		\newblock Evaluation of services using a fuzzy analytic hierarchy process.
		\newblock {\em Applied Soft Computing}, 5(1):23--33, 2004.
		
		\bibitem{mohtashami2021novel}
		Ali Mohtashami.
		\newblock A novel modified fuzzy best-worst multi-criteria decision-making
		method.
		\newblock {\em Expert Systems with Applications}, 181:115196, 2021.
		
		\bibitem{mou2016intuitionistic}
		Qiong Mou, Zeshui Xu, and Huchang Liao.
		\newblock An intuitionistic fuzzy multiplicative best-worst method for
		multi-criteria group decision making.
		\newblock {\em Information Sciences}, 374:224--239, 2016.
		
		\bibitem{rezaei2015best}
		Jafar Rezaei.
		\newblock Best-worst multi-criteria decision-making method.
		\newblock {\em Omega}, 53:49--57, 2015.
		
		\bibitem{rezaei2016best}
		Jafar Rezaei.
		\newblock Best-worst multi-criteria decision-making method: Some properties and
		a linear model.
		\newblock {\em Omega}, 64:126--130, 2016.
		
		\bibitem{rezaei2017multi}
		Jafar Rezaei, Alexander Hemmes, and Lori Tavasszy.
		\newblock Multi-criteria decision-making for complex bundling configurations in
		surface transportation of air freight.
		\newblock {\em Journal of Air Transport Management}, 61:95--105, 2017.
		
		\bibitem{roy1990outranking}
		Bernard Roy.
		\newblock The outranking approach and the foundations of electre methods.
		\newblock In {\em Readings in multiple criteria decision aid}, pages 155--183.
		Springer, 1990.
		
		\bibitem{royden1988real}
		Halsey~Lawrence Royden and Patrick Fitzpatrick.
		\newblock {\em Real analysis}, volume~32.
		\newblock Macmillan New York, 1988.
		
		\bibitem{saaty1990make}
		Thomas~L Saaty.
		\newblock How to make a decision: the analytic hierarchy process.
		\newblock {\em European journal of operational research}, 48(1):9--26, 1990.
		
		\bibitem{sevkli2010application}
		Mehmet Sevkli.
		\newblock An application of the fuzzy electre method for supplier selection.
		\newblock {\em International Journal of Production Research},
		48(12):3393--3405, 2010.
		
		\bibitem{vahidi2018sustainable}
		F~Vahidi, S~Ali Torabi, and MJ~Ramezankhani.
		\newblock Sustainable supplier selection and order allocation under operational
		and disruption risks.
		\newblock {\em Journal of Cleaner Production}, 174:1351--1365, 2018.
		
		\bibitem{van2017battle}
		Geerten Van~de Kaa, Daniel Scholten, Jafar Rezaei, and Christine Milchram.
		\newblock The battle between battery and fuel cell powered electric vehicles: A
		bwm approach.
		\newblock {\em Energies}, 10(11):1707, 2017.
		
		\bibitem{wan2021novel}
		Shuping Wan and Jiuying Dong.
		\newblock A novel extension of best-worst method with intuitionistic fuzzy
		reference comparisons.
		\newblock {\em IEEE Transactions on Fuzzy Systems}, 30(6):1698--1711, 2021.
		
		\bibitem{yang2019research}
		Qifeng Yang, Yingying Wang, and Yidong Ren.
		\newblock Research on financial risk management model of internet supply chain
		based on data science.
		\newblock {\em Cognitive Systems Research}, 56:50--55, 2019.
		
		\bibitem{ZADEH1965338}
		L.A. Zadeh.
		\newblock Fuzzy sets.
		\newblock {\em Information and Control}, 8(3):338--353, 1965.
		
		\bibitem{zekhnini2021analytic}
		Kamar Zekhnini, Anass Cherrafi, Imane Bouhaddou, and Youssef Benghabrit.
		\newblock Analytic hierarchy process (ahp) for supply chain 4.0 risks
		management.
		\newblock In {\em Artificial Intelligence and Industrial Applications: Smart
			Operation Management}, pages 89--102. Springer, 2021.
		
		\bibitem{zhao2018comprehensive}
		Haoran Zhao, Sen Guo, and Huiru Zhao.
		\newblock Comprehensive benefit evaluation of eco-industrial parks by employing
		the best-worst method based on circular economy and sustainability.
		\newblock {\em Environment, development and sustainability}, 20:1229--1253,
		2018.
		
		\bibitem{zheng2018smart}
		Pai Zheng, Honghui Wang, Zhiqian Sang, Ray~Y Zhong, Yongkui Liu, Chao Liu,
		Khamdi Mubarok, Shiqiang Yu, and Xun Xu.
		\newblock Smart manufacturing systems for industry 4.0: Conceptual framework,
		scenarios, and future perspectives.
		\newblock {\em Frontiers of Mechanical Engineering}, 13:137--150, 2018.
		
		\bibitem{zimmermann2011fuzzy}
		Hans-J{\"u}rgen Zimmermann.
		\newblock {\em Fuzzy set theory—and its applications}.
		\newblock Springer Science \& Business Media, 2011.
		
	\end{thebibliography}
\end{document}